\documentclass[review,3p]{elsarticle}
\usepackage{amssymb}

\usepackage{amsmath,amsfonts,bm}

\def\eqref#1{equation~\ref{#1}}

\def\1{\bm{1}}

\def\vb{{\bm{b}}}

\def\vf{{\bm{f}}}
\def\vg{{\bm{g}}}

\def\vn{{\bm{n}}}

\def\vu{{\bm{u}}}
\def\vv{{\bm{v}}}
\def\vw{{\bm{w}}}
\def\vx{{\bm{x}}}

\def\vz{{\bm{z}}}

\def\vH{\bm{H}}

\def\vL{\bm{L}}

\def\RR{\mathbb{R}}

\DeclareMathAlphabet{\mathsfit}{\encodingdefault}{\sfdefault}{m}{sl}
\SetMathAlphabet{\mathsfit}{bold}{\encodingdefault}{\sfdefault}{bx}{n}
\newcommand{\tens}[1]{\bm{\mathsfit{#1}}}

\newcommand{\R}{\mathbb{R}}

\newcommand{\bn}{{\bf n}} 

\newcommand{\bb}{{\bf b}}

\newcommand{\bx}{{\bf x}}

\newcommand{\cE}{{\cal E}}

\newcommand{\cM}{{\cal M}}

\newcommand{\cT}{{\cal T}}

\newcommand{\bsigma}{\mbox{\boldmath${\sigma}$}}

\newcommand{\bdelta}{\mbox{\boldmath${\delta}$}}
\newcommand{\bomega}{\mbox{\boldmath${\omega}$}}

\newcommand{\Om}{\Omega}

\newcommand{\C}{\rm I\kern-.5emC}
\newcommand{\G}{\Gamma}

\newcommand{\norm}[1]{\left\|{#1}\right\|}

\newcommand{\oversett}[2]{%
\mathop{#2}\limits^{\vbox to -.1ex{\kern -0.5ex\hbox{$\scriptstyle #1$}\vss}}}

\usepackage[utf8]{inputenc} 
\usepackage[T1]{fontenc}    
\usepackage{hyperref}       
\usepackage{url}            
\usepackage{booktabs}       
\usepackage{amsfonts}       
\usepackage{nicefrac}       
\usepackage{microtype}      
\usepackage{xcolor}         
\usepackage{amsmath, amsthm, amssymb}
\usepackage{graphicx}
\usepackage{subfigure}
\usepackage{wrapfig}
\usepackage{array}
\usepackage{multirow}
\usepackage{makecell}
\usepackage{algorithmic}
\usepackage{algorithm}

\newcolumntype{C}[1]{>{\centering\arraybackslash}p{#1}}

\newtheorem{prop}{Proposition}
\newtheorem{theo}{Theorem}
\newtheorem{lemm}{Lemma}

\begin{document}

\begin{frontmatter}
\title{Deep Ritz Method with Adaptive Quadrature for Linear Elasticity}
\author{Min Liu\fnref{mliu}}
\ead{liu66@purdue.edu}
\fntext[mliu]{School of Mechanical Engineering, Purdue University, 585 Purdue Mall,
West Lafayette, IN 47907-2088}
\author{Zhiqiang Cai\fnref{zcai}}
\ead{caiz@purdue.edu}
\fntext[zcai]{Department of Mathematics, Purdue University, 150 N. University Street, West Lafayette, IN 47907-2067}
\author{Karthik Ramani\fnref{mliu}}
\ead{ramani@purdue.edu}

\newcommand{\fix}{\marginpar{FIX}}
\newcommand{\new}{\marginpar{NEW}}

\begin{abstract}
In this paper, we study the deep Ritz method for solving the linear elasticity equation from a numerical analysis perspective. A modified Ritz formulation using the $H^{1/2}(\Gamma_D)$ norm is introduced and analyzed for linear elasticity equation in order to deal with the (essential) Dirichlet boundary condition. We show that the resulting deep Ritz method provides the best approximation among the set of deep neural network (DNN) functions with respect to the ``energy'' norm. \textcolor{black} {Furthermore, we demonstrate that the total error of the deep Ritz simulation is bounded by the sum of the network approximation error and the numerical integration error, disregarding the algebraic error.} To effectively control the numerical integration error, we propose an adaptive quadrature-based numerical integration technique with a residual-based local error indicator. This approach enables efficient approximation of the modified energy functional. Through numerical experiments involving smooth and singular problems, as well as problems with stress concentration, we validate the effectiveness and efficiency of the proposed deep Ritz method with adaptive quadrature. 
\end{abstract}



\end{frontmatter}
\section{Introduction}


In the past decade, Deep Neural Networks (DNNs) have achieved remarkable success in computer vision, natural language processing, and many other machine learning (ML) applications.  More recently, scientific machining learning methods based on DNN have also been applied
to modeling and solving complex engineering systems. These methods can be broadly divided into three categories based on how the DNN is used: 
(i) \textit{purely data-driven approaches}, which use supervised ML to create a surrogate model that regresses a physical model from a given simulation dataset or experimental observations \citep{Liang18,Gao20,Wang21,Badarinath21,iakovlev2021learning,li2021fourier}; (ii) \textit{physics-enhanced approaches}, which use semi-supervised ML to implement physical laws as a regularizing term  to solve a target regression problem with limited observation data \citep{raissi2017hidden,Long19,Khoo19} ; 
and (iii) \textit{ physics-driven approaches}, which impose physics into the loss functional and training process and rely on unsupervised ML to directly solve various types of PDEs \citep{Weinan17, Sirignano18,Berg18,Weinan18,Karniadakis19,CAI2020,Cai2021linear,LI-wei21}. 
Mathematically speaking, purely data-driven approaches can be thought of as methods for finding the best curve fit through the data points, using techniques such as least-squares regression. Physics-driven approaches, on the other hand, typically involve solving PDEs using numerical optimization methods. 
Numerous studies have demonstrated that DNNs possess highly desirable approximation properties that are not available in commonly used finite element methods. 
One such advantage is that a DNN can adapt its physical partition to match the underlying function being approximated \citep{LiuCai1,CAI2022,Cai2023}. In other words, a DNN model can dynamically adjust its function representation, akin to a moving mesh method, but without the necessity of a geometric mesh.

Due to the fact that the set of DNN functions does not form a linear space, existing physics-driven methods are typically based on either the energy minimization formulation \citep{Weinan18, Xu2020,LiuCai2} or various least-squares formulations \citep{Berg18, Karniadakis19, Sirignano18, CAI2020}. Energy minimization formulations are particularly suitable for problems that possess a natural minimization principle, such as many problems encountered in solid mechanics. On the other hand, the effectiveness of least-squares formulations depends on the specific principle employed. For example, the deep Galerkin method (DGM) \citep{Sirignano18} and the physics-informed neural networks (PINN) \citep{Karniadakis19} rely on the discrete $l^2$ norm least-squares principle, which applies to differential equations, and boundary and initial conditions. 
However, these methods suffer from suboptimal approximation and are limited to problems with $\vH^2$-smooth solutions, thereby excluding problems with geometric or interface singularities. To overcome these limitations, well-designed least-squares methods for PDEs can be employed, as described in books such as \cite{bg5} and papers such as \cite{Cai04, BeCaPa:19} for linear elasticity. Recently, the DNN-based first-order system least-squares (FOSLS) formulation has been utilized for the second-order elliptic PDEs \cite{CAI2020, HAGHIGHAT2021}.

In this paper, we aim to investigate the deep Ritz method for solving the linear elasticity equation from a numerical analysis perspective, as well as to introduce \textcolor{black}{a deep Ritz method with {\it adaptive} quadrature}. Originally proposed in \cite{Weinan18} for scalar elliptic PDEs, the deep Ritz method employs DNNs as the class of approximating functions and is based on the Ritz formulation of the underlying PDE. However, unlike finite element approximations, the deep Ritz method encounters two fundamental challenges arising from the characteristics of DNN functions. The first challenge pertains to enforcing the Dirichlet boundary condition effectively. The second challenge involves devising a numerical integration scheme that plays a crucial role in ensuring the accuracy and robustness of the DNN approximation to the solution of the underlying problem.

\textcolor{black}{Regarding the first issue, there are mainly two approaches. One is the Nitsche method \cite{Nitsche71, Liao21} 
that converts the Dirichlet boundary condition into the Robin boundary condition by {\it penalizing} the Neumann boundary condition, and the penalization constant has to be sufficiently small. This approach is equivalent to penalize the energy functional by the $L^2$ norm of the residual of the Dirichlet boundary condition \cite{Weinan18}. The other one is to enforce the Dirichlet boundary condition exactly through an auxiliary continuous function vanishing on the Dirichlet boundary\cite{URIARTE23}. In this paper, we explore the penalization method with $H^{1/2}$ norm that guarantees stability of the perturbed problem.}
Specifically, the standard minimization formulation is modified by adding the $\vH^{1/2}$ norm of the residual of the Dirichlet boundary condition to the energy functional.
To show the well-posedness of the modified minimization problem, we first prove a fundamental inequality of Korn's type in $\vH^1(\Omega)$ (see Lemma~1), and With this newly established inequality, 
we further show that the modified minimization problem has a unique solution and the solution continuously depends on the data (see Proposition~\ref{1}).
Based on the modified minimization problem, the deep Ritz method is defined by minimizing the modified energy functional over the set of DNN functions and the deep Ritz approximation with the exact integration and differentiation is proved to be the best approximation in the modified energy norm (see Theorem~1).

An evaluation of the modified energy functional includes integration over both the domain and the boundary, as well as differentiation at integration points. Naturally, the integration is approximated by quadrature-based methods, since the dimension of the linear elasticity problem is at most four (including space and time). Under a reasonable assumption on numerical integration, we demonstrate that the total error in the energy norm is bounded by the approximation error of the set of DNNs plus the numerical integration error (see Theorem~2). \textcolor{black}{It is important to note that solving the minimization problems under the deep Ritz formulation using DNNs gives rise to a high-dimensional and non-convex optimization problem. However, the algebraic error introduced during the solving/training process falls beyond the scope of this discussion.} 

In the finite element setting, it is trivial to control the numerical integration error because the unknown finite element approximation is a piece-wise polynomial on a {\it fixed} triangulation of the computational domain. However, controlling the numerical integration error for the deep Ritz method is difficult since the unknown DNN approximation is a composition function with several layers.  Moreover, the accuracy of the DNN approximation is determined by the quality of numerical integration mesh on which the solution can be approximated well by a selected quadrature rule \textcolor{black}{\cite{RIVERA22}.}
To overcome this obstacle, we propose an adaptive quadrature method that refines integration mesh and quadrature points adaptively. A modified residual-based local error indicator is used for marking subdomains to be refined.
The effectiveness and efficiency of the deep Ritz method with adaptive quadrature is studied for several benchmarks.

The rest of the paper is organized as follows. Section 2 introduces the modified Ritz formulation for linear elasticity problems and establishes its well-posedness (existence, uniqueness, and stability). Section 3 describes the discretization of the problem with DNN functions and shows the upper bounds of the approximation error. In Section 4, we propose the adaptive quadrature method and introduce the corresponding local error indicator. The last two sections present the numerical results and conclude the paper.

We will use the standard notation and definitions for the Sobolev space $\vH^{s}(\Omega)^d$ and $\vH^{s}(\Gamma)$ for a subset $\Gamma$ of the boundary of the domain $\Omega\in \RR^d$. The standard associated inner product and norms are denoted by $(\cdot,\cdot)_{s,\Omega,d}$ and $(\cdot,\cdot)_{s,\Gamma,d}$ and by $\Vert \cdot \Vert _{s,\Omega,d}$ and $\Vert \cdot \Vert _{s,\Gamma,d} $, respectively. When there is no ambiguity, the subscript $\Omega$ and $d$ in the designation of norms  will be suppressed. When $s=0$, $\vH^{0}(\Omega)^d$ coincides with $\vL^2(\Omega)^d$. In this case, the inner product and norm will be denoted by $(\cdot,\cdot)$ and $\|\cdot\|$.

\section{Modified Ritz Formulation of Linear Elasticity}
\label{formulation}
\textcolor{black}{Let $\Om$ be a bounded domain in $\RR^d$ ($d=2$ or $3$) with Lipschitz boundary $\partial\Omega={\Gamma}_D\cup{\Gamma}_N$, where $\Gamma_D$ and $\Gamma_N$ are disjoint.} 
Let $\vn$ be the outward unit
vector normal to the boundary. Denote by $\vu$ and $\tens{\sigma}$ the displacement field and the stress tensor.  
Consider the following linear elasticity problem
\begin{equation}
\label{pde}
\left\{
\begin{array}{rcll}
-\nabla \cdot \tens{\sigma} &=&\vf,  
\; &\text{ in }\, \Om,
\\ 
\tens{\sigma}(\vu) &=& 2\mu \tens{\epsilon}(\vu) + \lambda \nabla \cdot \vu \,\bdelta_{d\times d}\; &\text{ in }\, \Om
\end{array} 
\right.
\end{equation}
with boundary conditions $\vu \big|_{\G_D} =\vg_{_D}  
\quad\mbox{and }\,\, \big(\tens{\sigma}{\vn}\big)\big|_{\G_N} =\vg_{_N}$, 
where $\nabla \cdot$ is the divergence operator; $\tens{\epsilon}(\vu) = \frac{1}{2}\big(\nabla\vu + (\nabla\vu)^T\big)$ is the strain tensor;
the $\vf$, $\vg_{_D}$, and $\vg_{_N}$ are given vector-valued functions defined on
$\Om$, $\G_D$, and $\G_N$, representing body force, boundary displacement and boundary traction force condition respectively; $\bdelta_{d\times d}$ is the the $d$-dimensional identity matrix; $\mu$ and $\lambda$ are the material Lam$\acute{e}$ constants. 

Since it is difficult for directly constraining neural network functions to satisfy boundary conditions (see \cite{Weinan18}), as in \cite{CAI2020}, we enforce the Dirichlet (essential) boundary condition weakly using a half norm through the energy functional. 
The modified Ritz formulation of  problem~(\ref{pde}) is to find 
$\vu\in \vH^1(\Omega)^d$ such that 
\begin{equation}\label{mini}
    J(\vu)= \min_{\vv\in \vH^1(\Omega)^d} J(\vv),
\end{equation}
where the modified energy functional is given by
\begin{equation}\label{energy}
 J(\vv)=\dfrac12\left\{\int_{\Omega}\left(2\mu \left|\tens{\varepsilon}(\vv)\right|^2 + \lambda \left|\nabla\!\cdot\!\vv\right|^2\right)d\vx + \gamma\norm{\vv - \vg_{_D}}^2_{1/2,\Gamma_D} \right\} - (\vf,\,\vv) - (\vg_{_N}, \vv)_{0,{\Gamma_{N}}}.
\end{equation}  
Here, $\gamma=\mu\gamma_{_D}$ is a penalization constant scaled by $\mu$, and $\norm{\cdot}_{1/2,\Gamma_D}$ denotes the Sobolev-Slobodeckij norm given by
\begin{equation}\label{1/2norm}
\norm{\vv}_{1/2,\Gamma_D}=\left(\int_{\Gamma_D}|\vv|^2dS + \int_{\Gamma_D}\int_{\Gamma_D}\dfrac{|\vv(\vx)-\vv(\vx^\prime)|^2}{d(\vx,\vx^\prime)^{d}}dS(\vx)dS(\vx^\prime)\right)^{1/2},
\end{equation}
where $d(\vx,\vx^\prime)$ is the geodesic distance between $\vx$ and $\vx^\prime$ in $\Gamma_D$. Let 

\begin{eqnarray*}
a(\vu,\vv) 
=2\mu(\tens{\epsilon}(\vu),\tens{\epsilon}(\vv)) + \lambda (\nabla\!\cdot\! \vu,\nabla\!\cdot\! \vv) + \gamma\, (\vu,\vv)_{\frac12,\Gamma_D} \\ [2mm]
\,\mbox{ and }\, f(\vv)= (\vf,\vv)+(\vg_{_N},\vv)_{0,\Gamma_N} +  \gamma \,(\vg_{_D},\vv)_{\frac12,\Gamma_D},
\end{eqnarray*}

then the variational form of (\ref{mini}) is to finding $\vu\in \vH^1(\Omega)^d$ such that 
\begin{equation}\label{var}
     a(\vu,\vv)=f(\vv), \quad\forall\,\, \vv\in \vH^1(\Omega)^d.
\end{equation}

\begin{lemm}
For all $\vv\in \vH^1(\Omega)^d$, there exists a positive constant $C$ such that 
\begin{equation}\label{Korn}
    \norm{\vv}_{1,\Omega} \leq C\left(\|\tens{\varepsilon}(\vv)\|_{0,\Omega} + \|\vv\|_{1/2,\Gamma_D}
    \right).
\end{equation}
\end{lemm} 

\begin{proof}
For simplicity, we prove the validity of (\ref{Korn}) in $\RR^2$. To this end, denote the space of infinitesimal rigid motions in $\RR^2$ by $RM = \{\vv=(a,b)^T+c\,(y,x)^T |\, a, b, c \in \RR\}$.
For any $\vv\in \vH^1(\Omega)$, there exists a unique pair $(\vz,\vw)\in \hat{\vH}^1(\Omega)\times RM$ such that $\vv=\vz+\vw$,
where $\hat{\vH}^1(\Omega)=\left\{\vv\in \vH^1(\Omega)\big|\, \int_\Omega\vv={\bf 0},\,\, \int_\Omega\nabla\times\vv={\bf 0}
\right\}$. By the second Korn inequality and the fact that $\tens{\varepsilon}(\vw)={\bf 0}$ for any $\vw\in RM$, we have
\begin{equation}\label{2Korn}
    \norm{\vz}_{1,\Omega}\leq C\, \norm{\tens{\varepsilon}(\vz)}_{1,\Omega} =C\, \norm{\tens{\varepsilon}(\vv)}_{1,\Omega}.
\end{equation}
By the fact that $RM$ is a finite dimensional space, the triangle and trace inequalities give
\[
\norm{\vw}_{1,\Omega}\leq C\,\norm{\vw}_{1/2,\partial\Omega}\leq C\,\big(\norm{\vv}_{1/2,\partial\Omega}+\norm{\vz}_{1/2,\partial\Omega} \big) \leq C\,\big(\norm{\vv}_{1/2,\partial\Omega}+\norm{\vz}_{1,\Omega} \big).
\]
Now, (\ref{Korn}) is a direct consequence of the triangle inequality and the above two inequalities.
\end{proof}

\begin{prop}\label{1}
Problem {\em (\ref{mini})} has a unique solution $\vu\in \vH^1(\Omega)^d$. Moreover, the solution $\vu$ satisfies the following {\it a priori} estimate:
\begin{equation}\label{stable}
\|\vu\|_{1,\Omega}\leq C\, \left( \|\vf\|_{-1,\Omega} + \|\vg_{_D}\|_{1/2,\Gamma_D}
+ \|\vg_{_N}\|_{-1/2,\Gamma_N}\right).
\end{equation}
\end{prop}

\begin{proof}
\textcolor{black}{With the Korn inequality in (\ref{Korn}), the proof of the proposition is standard.}
\end{proof}

\section{Deep Ritz Neural Network Method}
\label{headings}

\begin{figure}[htbp]
\centering
\includegraphics[ height=1.8in]{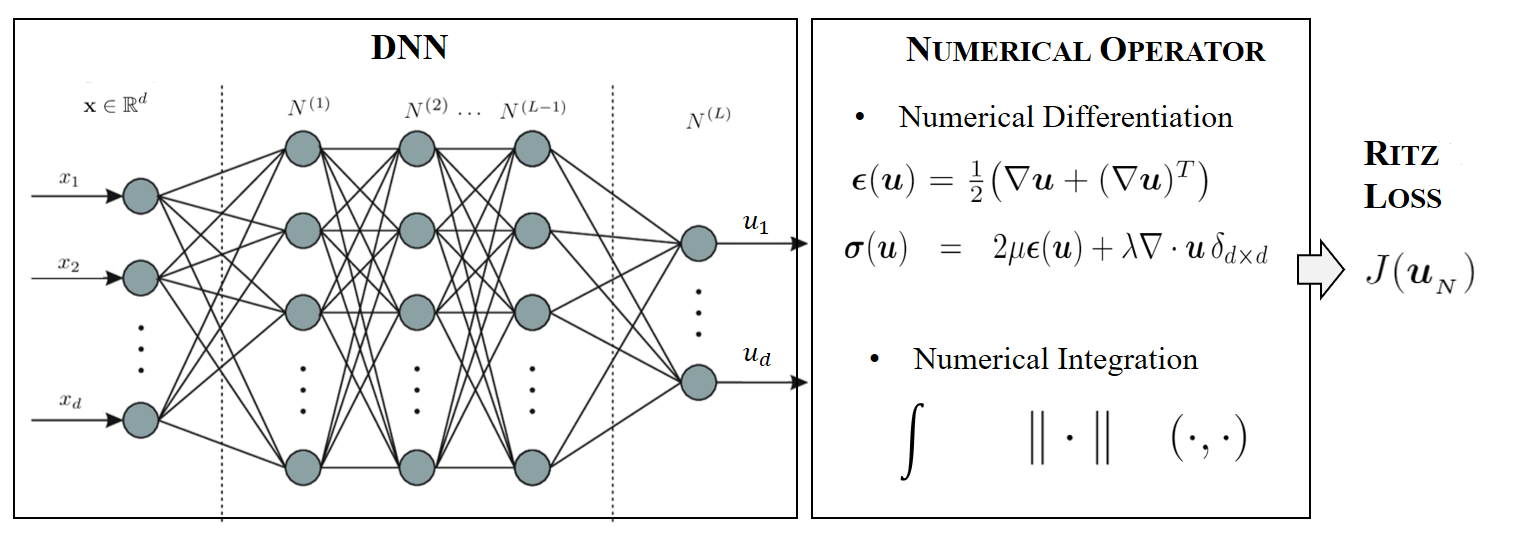}
\caption{Deep Ritz NN architecture. A fully connected $L$-layer network is employed to generate the map from an arbitrary spatial point $\vx$ in $\Om$ to its displacement $\vu(\vx)$, numerical operators are used to approximate the gradient, divergence and integral in the discrete energy functional $J(\vu_{_N})$ as the Ritz loss.}
\label{network1}
\end{figure}

\textcolor{black}{In this section, we describe the deep Ritz method which includes a standard fully connected DNN as the class of approximating functions and the discrete energy functional $J_{_\cT}(\vv)$ as an approximation of the energy functional $J(\vv)$ by numerical integration and differentiation.} The structure of the deep Ritz NN is illustrated in Figure \ref{network1}.

\subsection{Deep Neural Network}

For $j=1,\cdots,l-1$, let $N^{(j)}\!: \RR^{n_{j-1}}\! \rightarrow \RR^{n_{j}}$ be the vector-valued ridge function of the form
\begin{equation}\label{layerdef}
  N^{(j)}(\bx^{(j-1)})
  = \tau (\bomega^{(j)}\bx^{(j-1)}-\bb^{(j)})
  \quad\mbox{for }\,\, \bx^{(j-1)}\in \R^{n_{j-1}},
\end{equation}
where $\bomega^{(j)} \in \RR^{n_{j}\times n_{j-1}}$ and $\bb^{(j)}\in \R^{n_{j}}$ are the respective weights and bias to be determined; $\bx^{(0)}=\bx$; and $\tau(t)$ is a non-linear activation function. 
There are many activation functions  such as ReLU, ReLU$^p$, sigmoids, sinusoidal, and hyperbolic tangent. (see, e.g., \cite{pinkus1999}).

Let $\bomega^{(l)} \in \RR^{d\times n_{{l}-1}}$ and $\vb^{({l})} \in \RR^{d}$ be the output weights and bias. Then a $l$-layer neural network generates the following set of vector fields in $\RR^d$
\begin{equation}\label{NN}
{\cal {\scriptsize M}}_{N}(l)\!=\!\big\{
\bomega^{l}\!\left(\! N^{(l-1)}\! \circ \!\cdots \!\circ\!  N^{(1)}(\bx)\!\right) - \vb^{l}\! :\,  \bomega^{(j)} \!
\in \R^{n_{j}\times n_{j-1}}, \bb^{(j)}\in \R^{n_{j}} \mbox{ for all } j
\big\},
\end{equation}
where the symbol $\circ$ denotes the composition of functions. 

This class of approximating functions is rich enough to accurately approximate any continuous function defined on a compact set $\Omega\in \RR^d$ (see \cite{Cybenko1989, HornikS1989} for the universal approximation property). 
However, this is not the main reason why NNs are so effective in practice. 
One way to understand its approximation power is from the viewpoint of polynomial spline functions with free knots (\cite{Schumaker}). 
The set ${\cal {\scriptsize M}}_{N}(l)$ may be regarded as a beautiful extension of free knot splines from one dimensional scalar-valued function to multi-dimensional vector-valued function. It has been shown that the approximation of functions by splines can generally be dramatically improved if the knots are free.

\subsection{Discretization}

Note that neural network functions in ${\cal {\scriptsize M}}_{N}(l)$ are nonlinear with respect to the weights $\{\bomega^{(j)}\}_{j=1}^{l-1}$ and the bias $\{\bb^{(j)}\}_{j=1}^{l-1}$. This implies that it is difficult to discretize (\ref{pde}) by the conventional approach based on the corresponding variational formulation (\ref{var}). Instead, discretization using NNs should be based on an optimization formulation. In this paper, we employ the Ritz formulation (\ref{mini}) that minimizes the energy functional.

To approximate the solution of (\ref{pde}) using \textcolor{black}{a} neural network, the deep Ritz method \textcolor{black}{minimizes} the energy functional over the set ${\cal {\scriptsize M}}_{N}(l)$, i.e., 
finds $\vu_{_N}\in {\cM}_{_N}(l)\subset \vH^1(\Omega)^d$ such that 
\begin{equation}\label{Ritz}
J(\vu_{_N}) = \min_{\vv\in {\cM}_{_N}(l)} J(\vv).
\end{equation}
\textcolor{black}{Since $\cM^1_n(l)$ is not a linear space, problem (\ref{Ritz}) may have many solutions.  }

\begin{theo}
Let $\vu\in \vH^1(\Omega)^d$ be the solution of problem {\em (\ref{var})}, and let $\vu_{_N}\in {\cM}_{_N}(l)$ be a solution of {\em (\ref{Ritz})}. Then we have
 \begin{equation}\label{error-R}
     \|\vu-\vu_{_N}\|_a = \inf_{\vv\in {\cM}_{_N}(l)} \|\vu-\vv\|_a,
 \end{equation}
where $\|\vv\|_a:=\sqrt{a(\vv,\vv)}$ is the energy norm. 
\end{theo} 

\textcolor{black}{\begin{proof}
Since $\vu_{_N}\in {\cM}_{_N}(l)\subset \vH^1(\Omega)^d$, (\ref{error-R}) is a direct consequence of 
\[
\|\vu-\vu_{_N}\|_a^2=2\left(J(\vu_{_N})-J(\vu)\right)
\leq 2\left(J(\vv)-J(\vu)\right) = \|\vu-\vv\|_a^2 
\]
for any $\vv\in {\cM}_{_N}(l)$.
\end{proof}
}

Theorem~1 indicates that $\vu_{_N}$ is the best approximation with respect to the energy norm $\norm{\cdot}_a$, within the neural network functions class ${\cal {\scriptsize M}}_{N}(l)$ . 





\subsection{Numerical integration}
Evaluation of the energy functional requires integration and differentiation which are often computed numerically in practice.
This section discusses numerical integration schemes suitable for neural network functions. To this end, let us partition the domain $\Omega$ by a collection of subdomains
\[
{\cal T}=\{K\, :\, K\mbox{ is an open subdomain of } \Omega\}
\]
such that
\[
 \bar{\Omega} = \cup_{K\in {\cal T}} \bar{K}
 \quad\mbox{and}\quad
 K\cap T = \emptyset,
 \quad \forall\,\, K,\, T \in {\cal T}.
 \]
That is, the union of all subdomains of ${\cal T}$ equals to the whole domain $\Omega$, and any two distinct subdomains of ${\cal T}$ have no intersection. The resulting partitions of the boundary $\Gamma_{_D}$ and $\Gamma_{_N}$ are
\[
{\cal E}_{_D}=\{E
=\partial K \cap \Gamma_{_D}:\,\, K\in\mathcal{T}\}
\quad \mbox{and}\quad 
{\cal E}_{_N}=\{E
=\partial K \cap \Gamma_{_N}:\,\, K\in\mathcal{T}\}.
\]

When using a smooth activation function such as sigmoid, ReLU$^p$ etc., neural network functions belong to at least $C^1(\Omega)$, i.e., 
\begin{equation}\label{c1}
    {\cal {\scriptsize M}}_{N}(l) \subset C^1(\Omega).
\end{equation}
In this case, as in \cite{CAI2020} numerical integration may use the composite quadrature rule, such as composite mid-point, trapezoidal, Simpson, Gaussian, etc., defined on an {\it artificial} partition $\cT$ of the domain $\Omega$. Here the artificial partition refers its independence of the underlying geometry of the approximating function and it allows us to partition the domain with few restrictions. With a chosen numerical integration, differentiation is evaluated at quadrature points and can be done by either numerical differentiation with relatively small step size or automatic differentiation. 


\textcolor{black}{For simplicity of presentation, we describe the composite mid-point rule for interior and boundary integration in (\ref{energy}) and (\ref{1/2norm}) in two dimensions. Let $\vx_{_T}$ and $\bx_{_E}$ be the centroids of $T\in {\cal T}$ and $E\in {\cal E}_{_S}$ for $S=D$ and $N$. For any integrand $v(\vx)$, 
the composite ``mid-point'' quadrature rule over the domain $\Omega$ and the boundary $\Gamma_{_S}$ are given by
\[
\int_\Omega v(\vx) \,d\bx 
 \approx \sum_{T\in \cT} v(\vx_{_T}) \, |T|\quad\mbox{and}\quad
 \int_{\Gamma_S} v(\vx) \,ds 
 \approx \sum_{E\in {\cal E}_{_S}} v(\vx_{_E}) \, |E|,
\]
where $|T|$ and $|E|$ are the respective volume of element $T\in \cT$ and area of boundary element $E\in {\cal E}_{_S}$. 
}

\textcolor{black}{Let the Dirichlet boundary $\Gamma_D$ be the union of several disjoint $\Gamma_D^k$ for $k=1,\ldots,I$. On each $\Gamma^k_D$, to approximate the Sobolev-Slobodeckij semi-norm in two dimensions, denote by $\vv[\vx;\vx^\prime]$ the integrand in the second term of (\ref{1/2norm}) as the divided difference of a vector-valued function $\vv(\vx)$ along $\Gamma^k_D$, where each component of $\vv[\vx;\vx^\prime]$ is given by
\[
v_i[\vx;\vx^\prime]=\left\{\begin{array}{ll}
 \dfrac{d v_i(\vx)}{d \Gamma^k_D},   &  \vx^\prime=\vx,\\ [4mm]
 \dfrac{v_i(\vx)-v_i(\vx^\prime)}{d(\vx,\vx^\prime)},    &  \vx^\prime\not=\vx.
\end{array}\right. 
\]
Then we have 
\begin{equation}\label{1/2-d}
 \int_{\Gamma^k_D}\int_{\Gamma^k_D} \dfrac{|\vv(\vx)-\vv(\vx^\prime)|^2}{d(\vx,\vx^\prime)^2}ds(\vx)ds(\vx^\prime) \,
 \approx \sum_{E\in {\cal E}^k_{_D}} \sum_{E^\prime\in {\cal E}^k_{_D}}|\vv[\vx_{_E};\vx_{_{E^\prime}}]|^2 \, |E|\,|E^\prime|.
\end{equation}
 }

\textcolor{black}{
Define the discrete bilinear and linear forms, $a_{_\cT}(\cdot,\cdot)$ and $f_{_\cT}(\cdot)$, by
\begin{eqnarray*}
 a_{_\cT}(\vu,\!\vv)
&=& 2\mu\sum_{T\in\cT} \big(\tens{\varepsilon}(\vu) : \tens{\varepsilon}(\vv)\big) (\vx_{_T})  + \lambda \sum_{T\in\cT}\big(\nabla\!\cdot\!\vu \nabla\cdot\vv\big)(\vx_{_T})\qquad\\[2mm]
&&\quad  +\mu{\gamma}_{_D}  \left\{\sum_{E\in\cE_{_D}}(\vu\cdot\vv)(\vx_{_E})
+ \sum_{k=1}^I\sum_{E\in {\cal E}^k_{_D}} \sum_{E^\prime\in {\cal E}^k_{_D}}\vu[\vx_{_E};\vx_{_{E^\prime}}]\cdot\vv[\vx_{_E};\vx_{_{E^\prime}}]  |E||E^\prime|\right\}\qquad\qquad \qquad\\[2mm]
f_{_\cT}(\vv)
&=& \sum_{T\in\cT}(\vf\!\cdot\!\vv)(\vx_{_T}\!) +
\sum_{E\in\cE_{_N}}(\vg_{_N}\!\cdot\!\vv)(\vx_{_E}) 
\qquad\\[2mm]
&&\quad  + \mu\gamma_{_D} \left\{
\sum_{E\in\cE_{_D}}(\vg_{_D}\!\cdot\!\vv)(\vx_{_E}) + \sum_{k=1}^I
\sum_{E\in {\cal E}^k_{_D}} \sum_{E^\prime\in {\cal E}^k_{_D}}\vg_{_D}[\vx_{_E};\vx_{_{E^\prime}}]\!\cdot\!\vv[\vx_{_E};\vx_{_{E^\prime}}] |E||E^\prime|\right\}.
\end{eqnarray*}
}
Define the discrete counterpart of the energy function $J(\cdot)$ by 
\begin{eqnarray}\nonumber
J_{_\cT}(\vv)=\dfrac12 a_{_\cT}(\vv,\vv) - f_{_\cT}(\vv).
\end{eqnarray}
Then the deep Ritz approximation to the solution of (\ref{pde}) is to seek $\vu_{_\cT}\in {\cM}_{_N}(l)$ such that 
\begin{equation}\label{Ritz-n-Q}
J_{_\cT}(\vu_{_\cT}) = \min_{\vv\in {\cM}_{_N}(l)} J_{_\cT}(\vv).
\end{equation}
To understand the effect of numerical integration, we extend the first Strang lemma for the Galerkin approximation over a subspace (see, e.g, \cite{Ciarlet78}) to the Ritz approximation over a subset. 



\textcolor{black}{
\begin{theo}\label{t:Q}
Assume that there exists a positive constant $\beta$ independent of ${\cM}_{_N}(l)$ such that
 \begin{equation}\label{Q}
     \beta\, \norm{\vv}_a^2 \leq a_{_\cT}(\vv,\vv), \quad\forall\,\, \vv\in {\cM}_{_N}(l).
 \end{equation}
Let $\vu$ be the solution of {\em (\ref{mini})} and $\vu_{_\cT}$ a solution of {\em (\ref{Ritz-n-Q})}. Then we have
\begin{equation}\label{E-Q}  
 \norm{\vu-\vu_{_\cT}}_a
\leq \! \dfrac{2}{\beta}\!\! \sup _{\vw\in {\cM}_{2N}(l)} \!\!\dfrac{|f(\vw) - f_{_\cT}(\vw)|}{\|\vw\|_a}
+ \dfrac{2+\beta}{\beta}\!\!
\inf_{\vv\in {\cM}_{_N}(l)} \left\{\norm{\vu-\vv}_a + \!\!\sup _{\vw\in {\cM}_{2N}(l)}\!\! \dfrac{|a(\vv,\vw)-a_{_\cT}(\vv,\vw)|}{\|\vw\|_a}\right\}.  
\end{equation}
\end{theo}
\begin{proof}
For any $\vv\in {\cM}_{_{N}}(l)\subset \vH^1(\Omega)^d$, we have
\[
 J_{_\cT}(\vu_{_\cT}) \leq J_{_\cT}(\vv)
 \quad\mbox{and}\quad
 a(\vu,\vu_{_\cT}-\vv)=f(\vu_{_\cT}-\vv).
 \] 
It follows from assumption (\ref{Q}) and the definition of $J_{_\cT}(\cdot)$ that
\begin{eqnarray*}
&& \dfrac{\beta}{2} \norm{\vu_{_\cT}-\vv}^2_a
\leq  \dfrac12 a_{_\cT}(\vu_{_\cT}-\vv,\vu_{_\cT}-\vv) \\[2mm] \nonumber
&=& J_{_\cT}(\vu_{_\cT}) -J_{_\cT}(\vv) +f_{_\cT}(\vu_{_\cT}-\vv) -a_{_\cT}(\vv,\vu_{_\cT}-\vv)  
\leq  f_{_\cT}(\vu_{_\cT}-\vv)-a_{_\cT}(\vv,\vu_{_\cT}-\vv)  \\[2mm] \nonumber
&=& \Big(f_{_\cT}(\vu_{_\cT}-\vv) -f(\vu_{_\cT}-\vv)\Big) +\Big(a(\vv,\vu_{_\cT}-\vv)-a_{_\cT}(\vv,\vu_{_\cT}-\vv)\Big)+ a(\vu-\vv, \vu_{_\cT}-\vv).
\end{eqnarray*}
which, together with the triangle and Cauchy-Schwarz inequalities, implies
\begin{eqnarray*}
\dfrac{\beta}{2} \norm{\vu_{_\cT}-\vv}_a
&\leq &\dfrac{|f_{_\cT}(\vu_{_\cT}-\vv) -f(\vu_{_\cT}-\vv)|}{\norm{\vu_{_\cT}-\vv}_a} +\dfrac{|a(\vv,\vu_{_\cT}-\vv)-a_{_\cT}(\vv,\vu_{_\cT}-\vv)|}{\norm{\vu_{_\cT}-\vv}_a} + \norm{\vu-\vv}_a  \\[4mm]   
&\leq & \sup _{\vw\in {\cM}_{{2N}}(l)}\dfrac{|f_{_\cT}(\vw) -f(\vw)|}{\norm{\vw}_a} +\sup _{\vw\in {\cM}_{{2N}}(l)}\dfrac{|a(\vv,\vw)-a_{_\cT}(\vv,\vw)|}{\norm{\vw}_a} + \norm{\vu-\vv}_a.
\end{eqnarray*}
Combining the above inequality with the triangle inequality
 \[
 \norm{\vu-\vu_{_\cT}}_a \leq \norm{\vu-\vv}_a + \|\vv-\vu_{_\cT}\|_a
 \]
and taking the infimum over all $\vv\in {\cM}_{_N}(l)$ yield (\ref{E-Q}). This completes the proof of the theorem.
\end{proof}
}
This theorem indicates that the total error in the energy norm is bounded by the approximation error of the set of neural network functions plus the numerical integration error. 

\section{Adaptive Quadrature Method}

As indicated in Theorem~\ref{t:Q}, numerical integration plays an important role in NN-based numerical methods. How to control the numerical integration error for the deep Ritz method is a non-trivial matter because the unknown DNN approximation is a composition function with several layers. To overcome this obstacle, in this section we propose adaptive Ritz method that refines integration mesh adaptively.  

Numerical integration defined in the previous section is based on an artificial partition $\cT$ of the domain $\Omega$. This partition may not capture well the variation of the underlying solution and hence (\ref{Ritz-n-Q}) would possibly lead to an inaccurate approximation. 

One may choose a uniform partition with sufficient fine mesh; however, it is cost inefficient. In this section, we describe an adaptive quadrature algorithm on numerical integration introduced in \cite{LiuCai1} under the assumption that the neural network is large enough to approximate the solution accurately.

A key ingredient for an adaptive quadrature scheme is an efficient local error indicator. In this paper, we use a modified residual-based indicator. To this end, let $\cT$ be the current integration mesh and $\vu_{_\cT}$ be a solution of (\ref{Ritz-n-Q}). For each $T\in\cT$, we define the following local error indicator for each $T\in\cT$,
\begin{equation}\label{indicator}
 \eta_{_T}(\vu_{_\cT})=\big|T\big|^{1/d}\big\|\nabla \cdot \bsigma_{_\cT} +\vf\big\|_{0,T},
\end{equation}
where $\bsigma_{_\cT}$ is the numerical stress given by
\[
\bsigma_{_\cT}=2\mu \tens{\epsilon}\big(\vu_{_\cT}\big) + \lambda \nabla \cdot \vu_{_\cT} \,\bdelta_{d\times d}.
\]
Note that the typical jump terms in finite element vanish due to the fact that ${\cM}_{_N}(l)$ is in $C^1(\Omega)$. The $L^2$ norm of the residual $\nabla \cdot \bsigma_{_\cT} +\vf$ on each $T\in\cT$ may be approximated as follows,
\[
\big\|\nabla \cdot \bsigma_{_\cT} +\vf\big\|_{0,T}\approx |T|^{-\frac12}\left|\int_T\big(\nabla \cdot \bsigma_{_\cT} +\vf\big)\,d\vx\right| = |T|^{-\frac12}\left|\int_{\partial T}\bsigma_{_\cT}\bn\,dS + \int_T\vf\,d\vx\right|
\]
which implies
\begin{equation}\label{indicator2}
 \eta_{_T}(\vu_{_\cT})\approx \big|T\big|^{\frac{2-d}{2d}}\left|\int_{\partial T}\bsigma_{_\cT}\bn\,dS + \int_T\vf\,d\vx\right|.
\end{equation}


With this local error indicator, we then define a subset $\hat{\cT}$ of $\cT$ by using either the following bulk marking strategy: finding a minimal subset $\hat{\cT}$ of $\cT$ such that
\begin{equation}\label{marking-1}
    \sum_{T\in \hat{\cT}} \eta^2_{_T}
    \ge \gamma_1\, \sum_{T\in \cT} \eta^2_{_T}
    \quad\mbox{for }\,\, \gamma_1\in (0,\,1)
\end{equation}
or the average marking strategy:
 \begin{equation}\label{marking-2}
    \hat{\cT} =\left\{T\in \cT\, :\,
    \eta_{_T}
    \ge \, \dfrac{\gamma_2}{\#\cT}\sum_{T\in \cT}\eta_{_T}\right\}, \quad\mbox{for }\,\, \gamma_2\in \left(0,\,\dfrac{\max\{\eta_{_T}\}}{\#\cT}\sum_{T\in \cT}\eta_{_T}\right)
\end{equation}
where $\#\cT$ is the number of subdomains of $\cT$. For each marked domain $T\in\hat{\cT}$, we subdivide $T$ into $2^d$ subdomains uniformly and denote the refinement partition by $\cT^\prime$. 

Let $\vu_{_\cT}$ be an approximation of (\ref{Ritz-n-Q}) based on an initial partition $\cT$, which in general, is an uniform partition of the domain, the adaptive quadrature refinement method is summarized as follows,

\begin{algorithm}
{\bf {\sc \bf Algorithm 3.1}} Adaptive Quadrature Refinement (AQR) with a fixed NN.
\begin{itemize}
     \item[(1)] for each $T\in\cT$, compute the local error indicator $\eta_{_\cT}$;
    \item[(2)] mark $\cT$ by the marking strategy and refine marked element to obtain a new partition $\cT^\prime$;
    \item[(3)] numerically solve the minimization problem in (\ref{Ritz-n-Q}) on $\cT^\prime$;
    \item[(4)] if $\eta (\vu_{_{\cT^\prime}}) \leq \gamma \eta (\vu_{_{\cT}})$,
    go to Step (1) with $\cT = \cT^\prime$; otherwise, output $\cT$.
\end{itemize}
\end{algorithm}

\smallskip

As indicated in \cite{LiuCai1}, the stopping criterion used in Algorithm~3.1 is based on whether or not the quadrature refinement on numerical integration improves approximation accuracy.  When the refinement does not improve accuracy much, the adaptive quadrature stops and outputs the current integration mesh.

\section{Numerical Studies}
\label{experiments}
In this section, we present our numerical results for several 2D problems. In all experiments, the DNN structure is \textcolor{black}{represented} as $d_{in}$-$n_1$-$n_2\cdots n_{l-1}$-$d_{out}$ for a $l$-layer network with $n_1$, $n_2$ and $n_{l-1}$ neurons in the respective first, second, and $(l-1)$th hidden layers, and $d_{in}=d_{out}$ represent the network input and output dimensions. The minimization of the deep Ritz NN loss function (\ref{Ritz-n-Q}) is solved using the Adam version of gradient descent \citep{kingma2015}. All differential operators are calculated using numerical differentiation (ND) with the step size $\Delta x = h/4$, where $h= \min\{\big|T\big|^{1/d}\}$ is the smallest partition size of an adaptive integration mesh. \textcolor{black}{All experiments use sigmoid ($\sigma (t) = \frac{1}{1+e^{-t}}$) as the activation function. For adaptive quadrature, the average marking strategy is reported due to its computational simplicity.}

\subsection{Test case I: smooth stress distribution}
Consider problem (\ref{pde}) defined on $\Omega=(-1,1)\times (-1,1)$ with the body force
\[
\vf=2\mu \big(3-x^2-2y^2-2xy,3-2x^2-y^2-2xy\big)^T + 2\lambda \big(1-y^2-2xy,1-x^2-2xy\big)^T,
\]
and the traction 
\[
\vg_{_N}=2(y^2-1)\big(2\mu+\lambda, \mu \big)^T
\]
on $\Gamma_{_N}=\{(1,y):\, y\in (-1,1)\}$, with the clamped boundary condition on $\Gamma_{_D}=\partial\Omega\setminus \Gamma_{_N}$. The exact solution of the test problem has the form
\[
\vu(x,y)=(1-x^2)(1-y^2)\big(1,1\big)^T.
\]



\begin{figure}[htb!]
  \centering 
\subfigure[The Problem]{ 
    \label{test1ritz:a} 
    \includegraphics[width=1.5in]{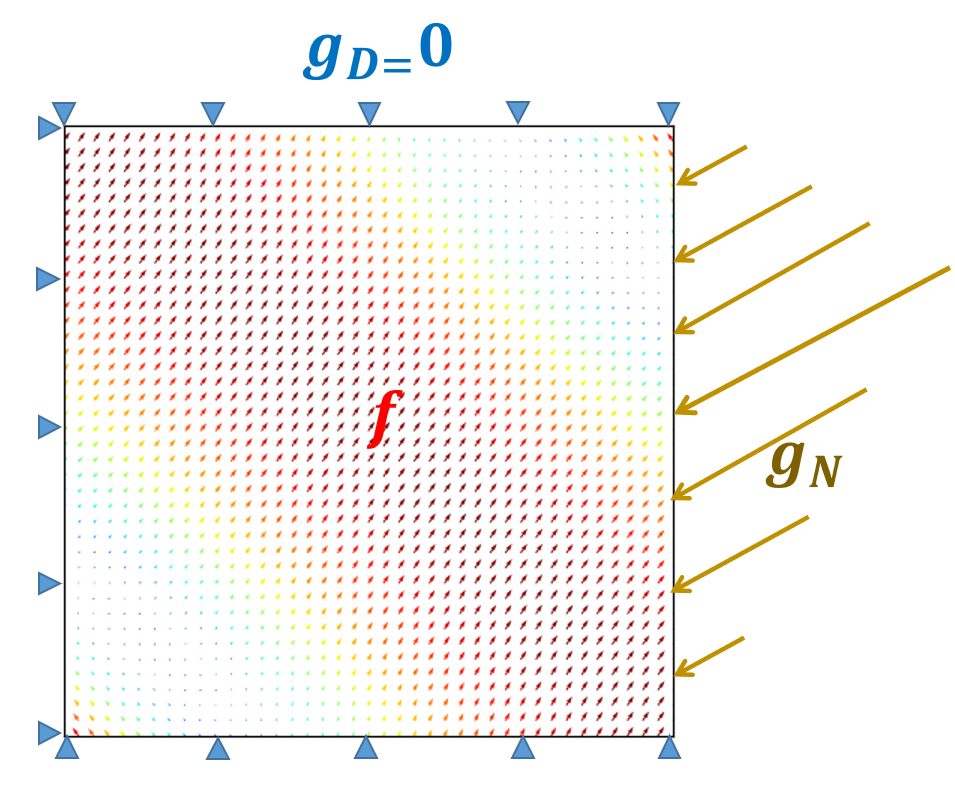}}
  \subfigure[Displacement ${u_x}$]{ 
    \label{test1ritz:b} 
    \includegraphics[width=1.5in]{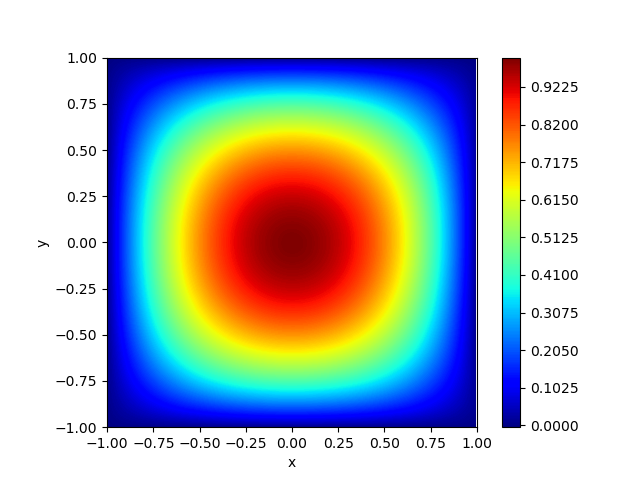}}
    \hspace{0.1in}
   \subfigure[Displacement ${u_y}$]{ 
    \label{test1ritz:c} 
    \includegraphics[width=1.5in]{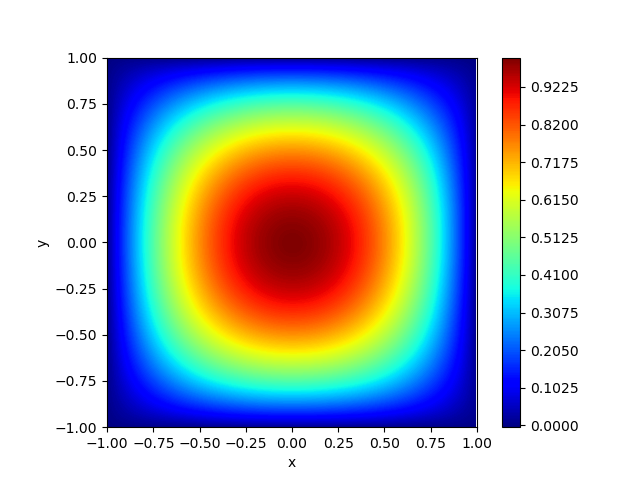}}
    \hspace{0.1in}
  \subfigure[Stress ${\sigma}_{xx}$]{ 
    \label{test1ritz:d} 
    \includegraphics[width=1.5in]{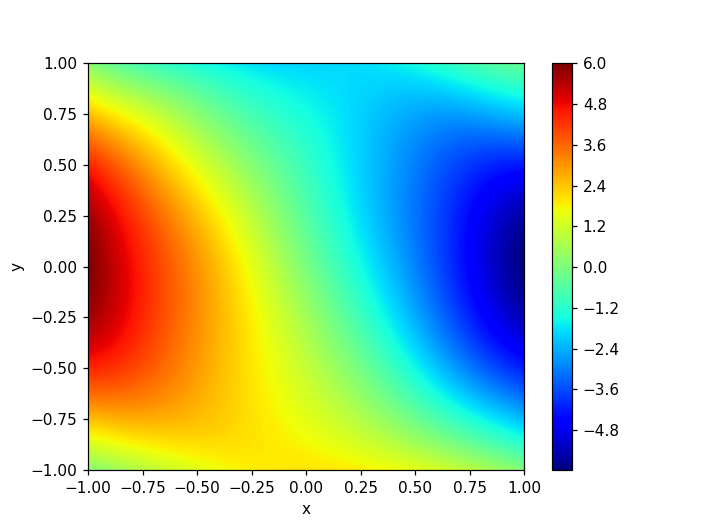}}
    \hspace{0.1in}
  \subfigure[Stress ${\sigma}_{yy}$]{ 
    \label{test1ritz:e} 
    \includegraphics[width=1.5in]{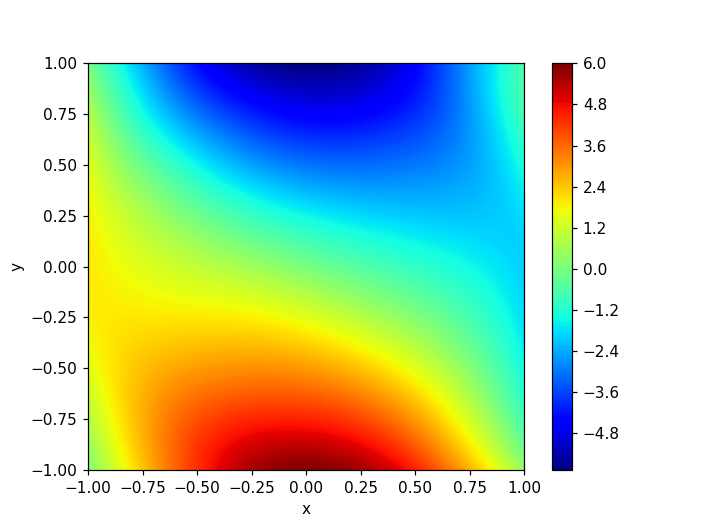}}
     \hspace{0.1in}
  \subfigure[Stress ${\sigma}_{xy}$]{ 
    \label{test1ritz:f} 
    \includegraphics[width=1.5in]{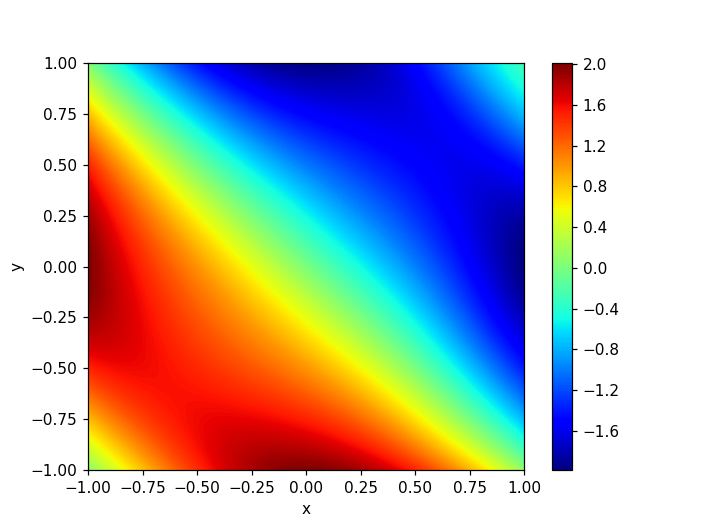}}

  \caption{Test case I problem ($\mu=1$, and $\lambda=1$) and numerical results using deep Ritz (2-32-32-2) with fixed $400\times400$ quadrature points.) }
  \label{Result-1ritz} 
\end{figure}

Set the material property $\mu=1$, and $\lambda=1$, we first test three-layer DNNs of varying number of neurons and different numerical quadrature resolutions. Uniformly distributed quadrature points of size $100\times100$, $200\times200$ and $400\times400$ are used to evaluate the effect of numerical integration combined with three network structures. 
Table. \ref{test1_table1} list the numerical results. As shown in the table, With a small DNN of 106 parameters (2-8-8-2) and $100\times100$ uniformly distributed quadrature points, deep Ritz can approximate the problem at a relative energy norm of 0.1658. Increasing the resolution of quadrature reduces the numerical integration error and therefore improves the approximation accuracy. E.g. by increasing the number of quadrature points to $200\times200$, and further to $400\times400$, the accuracy of the numerical solution is continuously improved, as measured by the relative energy norm, or the $L2$ norm of $\vu$ and $\sigma$.  On the other hand, larger DNNs with more parameters have better expressive power and thus can approximate the solution with better accuracy, this is a property from Theorem 2 and is  confirmed experimentally as the results given in Table. \ref{test1_table1}. A three-layer DNN (2-32-32-2) combined with finer quadrature points $400\times400$ has a better accuracy compared with smaller network or coarser quadrature points. The result is also depicted graphically in  Figure. \ref{Result-1ritz}. 
 


\begin{table}[htbp]
\centering
\caption{Numerical results of deep Ritz method for test case I using fixed quadrature}
\vspace{5pt}
\begin{tabular}{|*{5}{C{.8in}|}}
\hline
{\bfseries \makecell{DNN \\(No. para.)} } & {\bfseries $\#\cT$}
   &  $\dfrac{\|\vu-\vu_{_N}\|_a}{\|\vu\|_a}$ &  $\dfrac{\|\sigma-\sigma_{_N}\|}{\|\sigma\|}$ &  $\dfrac{\|\vu-\vu_{_N}\|}{\|\vu\|}$ \\
\hline
\multirow{3}{*}{\makecell{2-8-8-2 \\(106)}} & $100\times100$ & 16.58\% & 16.35\%  & 6.42\%  \\
\cline{2-5}
 &  $200\times200$ & 10.81\% & 10.61\%   & 3.91\% \\
\cline{2-5}
 &  $400\times400$ & 6.09\% & 6.00\%   & 2.15\% \\
\hline

\multirow{3}{*}{\makecell{2-16-16-2 \\(338)}} & $100\times100$ & 11.94\%  & 11.77\%    & 4.17\%  \\
\cline{2-5}
 &  $200\times200$ & 8.55\%  & 8.50\%    & 2.90\% \\
\cline{2-5}
 &  $400\times400$ & 5.94\%  & 5.93\%   & 1.96\%  \\
\hline

\multirow{3}{*}{\makecell{2-32-32-2 \\(1186)}} & $100\times100$ & 3.3\%  & 3.67\%    & 1.42\%  \\
\cline{2-5}
 &  $200\times200$ & 2.76\%  & 2.73\%   & 1.13\% \\
\cline{2-5}
 &  $400\times400$ & 2.19\% & 2.17\%  & 0.93\% \\
\hline
\multicolumn{5}{l}{\small *Training details: $\gamma_{_D}= 100$;} \\
\multicolumn{5}{l}{\small 
total number of iterations for each row: 200,000;} \\
\multicolumn{5}{l}{\small learning rate initial is 0.01 and it decays $90\%$ every 50000 iterations.}\\
\end{tabular}
\label{test1_table1}
\end{table}

Second, we test the adaptive quadrature refinement (AQR) method using the DNN structure 2-32-32-2. Starting from the initial $100\times100$ uniformly distributed quadrature points, and using the residual-based local error indicator(\ref{indicator}) and the average marking strategy(\ref{marking-2}) with $\gamma_2=1$, the AQR process stops at run 4 and reaches a relative energy norm 0.0174. As shown in Table \ref{test1_table2}, with 45,757 quadrature points after three-run AQR, the adaptive deep Ritz can achieve a similar relative energy norm as the fixed uniform quadrature method using 160,000 quadrature points. Adding more runs of adaptive quadrature process will improve the approximation accuracy, but it converges to a limit when the network approximation error becomes dominant. For example, by increasing the number of quadrature points from 45,757 to 99,262, error measured by the three norms are not improving significantly. To further improve the accuracy, one may need to enlarge the DNN size to obtain a better network approximation power.

\begin{table}[htb]
\centering
\caption{Numerical results of deep Ritz for experiment 1 using AQR.}
\vspace{5pt}
\begin{tabular}{|c|c|c|c|c|c|}
\hline
{\bfseries AQR run} & {\bfseries $\#\cT$}
   &  $\dfrac{\|\vu-\vu_{_N}\|_a}{\|\vu\|_a}$ &  $\dfrac{\|\sigma-\sigma_{_N}\|}{\|\sigma\|}$ &  $\dfrac{\|\vu-\vu_{_N}\|}{\|\vu\|}$  & $\dfrac{\sum_{T\in \cT}\eta_{_T}}{\#\cT}$\\
\hline
1 & $10,000$ & 3.73\% & 3.67\%  & 1.42\% & 0.0004\\
\hline
2 & $21,145$ & 2.81\% & 2.78\%   & 1.17\% & 0.0002\\
\hline
3 & $45,757$ & 1.98\% & 1.97\%   & 0.89\% &5e-5\\
\hline
4 & $99,262$ & 1.74\% & 1.73\%   & 0.82\% &2e-5 \\
\hline
\multicolumn{6}{l}{\small *Training details: $\gamma_{_D}= 100$;} \\
\multicolumn{6}{l}{\small run $2-4$ are trained using weight transferred from the previous run; } \\
\multicolumn{6}{l}{\small each trained 100,000 iterations using fixed learning rate 0.001.}\\
\end{tabular}
\centering
\label{test1_table2}
\end{table}

\subsection{Test case II: L-shape plate with corner singularity}
The second test is a common benchmark problem with a re-entrant corner forming a typical point singularity \citep{Harper19}. The problem is posed on an $L-$shaped domain $\Omega = (-1,1)^2\setminus([0,1]\times[-1,0])$ with a body force $\vf={\bf 0}$. The known analytical solution is,
\[
\vu = [A\cos\theta-B\sin\theta, A\sin\theta+B\cos\theta]^T,
\]
where $r, \theta$ are the polar coordinates and 
\begin{equation*}
\left\{
\begin{array}{l}
A=\frac{r^\alpha}{2\mu}\Big(-(1+\alpha)\cos\big((1+\alpha)\theta\big) + C_1(C_2-1-\alpha)\cos\big((1-\alpha)\theta\big)\Big),
\\[2mm]
B=\frac{r^\alpha}{2\mu}\Big((1+\alpha)\sin\big((1+\alpha)\theta\big) - C_1(C_2-1+\alpha)\sin\big((1-\alpha)\theta\big)\Big).
\end{array} 
\right.
\end{equation*}
Here the critical exponent $\alpha\approx0.544483737$ is the solution of the equation $\alpha\sin{(2\omega)} + \sin{(a\omega\alpha)} =0$ with $\omega=3\pi/4$ and $C_1=-(\cos{((\alpha+1)\omega})/(\cos{((\alpha-1)\omega})$, $C_2=(2(\lambda+2\mu))/(\lambda+\mu)$ \citep{Alberty02}. A bronze  material with Young's modulus $E =100000$ and Poisson's ratio $\nu=0.3$ is tested with the Neumann BCs prescribed on $\Gamma_{_N}=\{(1,y):\, y\in (0,1)\}$ and Dirichlet BCs on $\Gamma_{_D}=\partial\Omega\setminus \Gamma_{_N}$.

\begin{figure}[htb!]
  \centering 
\subfigure[Displacement ${u_x}$]{ 
    \label{test2:u1} 
    \includegraphics[width=1.5in]{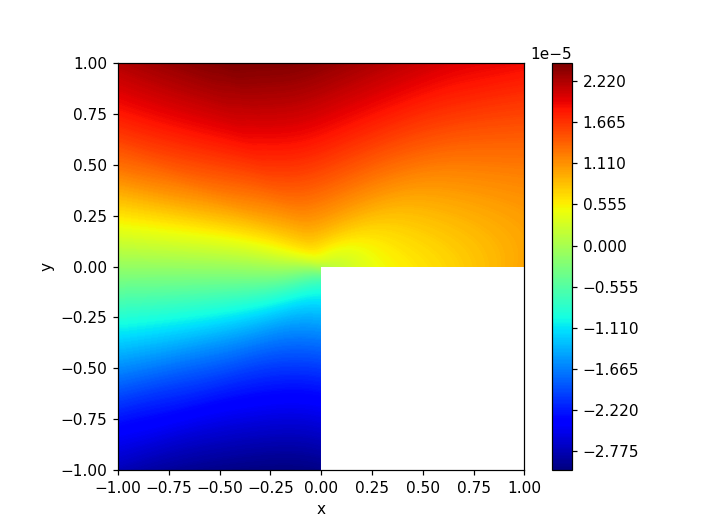}}
  \subfigure[Displacement ${u_y}$]{ 
    \label{test2:u2} 
    \includegraphics[width=1.5in]{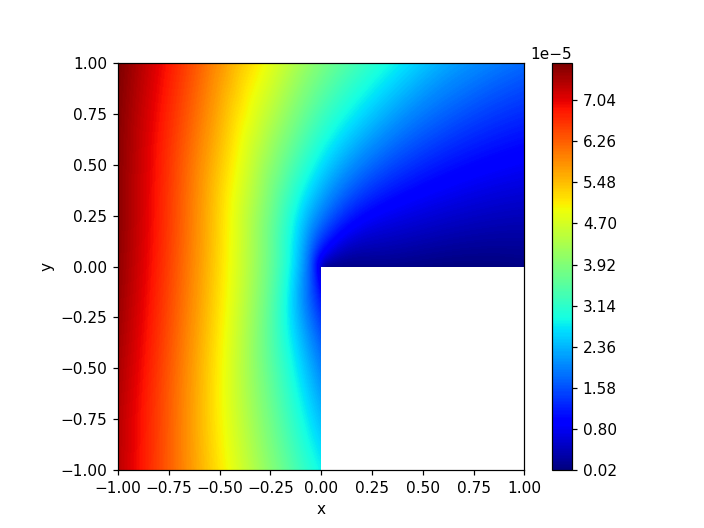}}
    \hspace{0.1in}
   \subfigure[Stress ${\sigma}_{xx}$]{ 
    \label{test2:s11} 
    \includegraphics[width=1.5in]{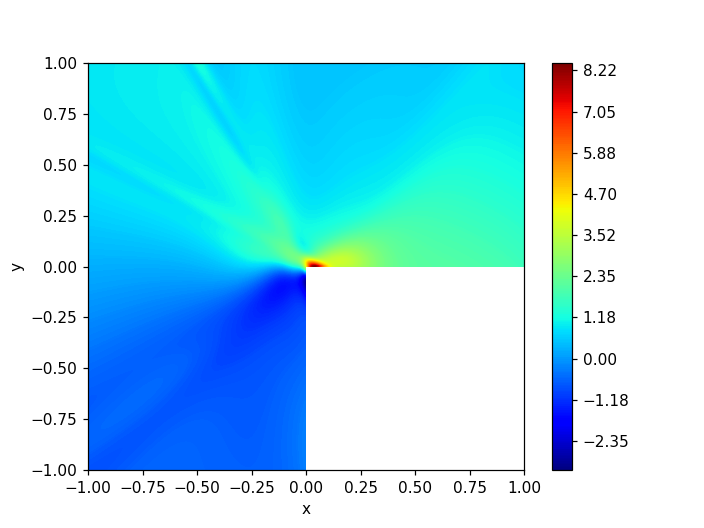}}
    \hspace{0.1in}
  \subfigure[Stress ${\sigma}_{yy}$]{ 
    \label{test2:s22} 
    \includegraphics[width=1.5in]{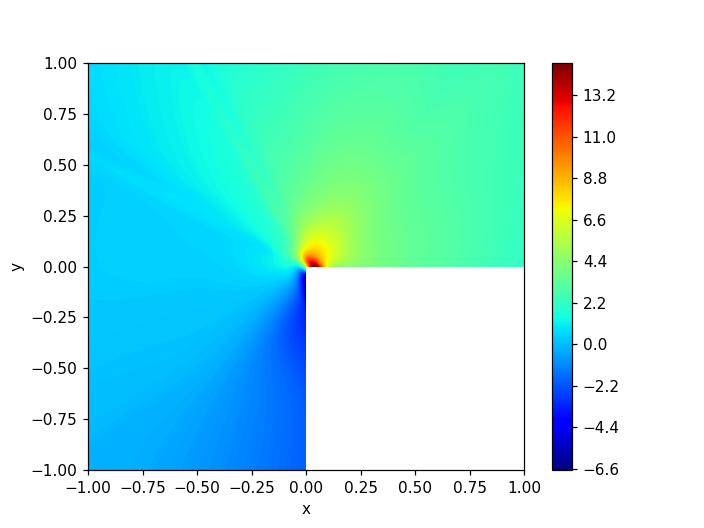}}
    \hspace{0.1in}
  \subfigure[Stress ${\sigma}_{xy}$]{ 
    \label{test2:s12} 
    \includegraphics[width=1.5in]{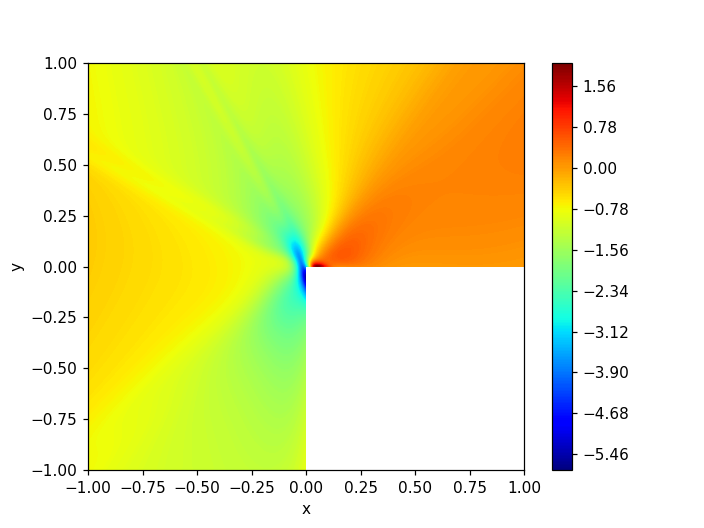}}
     \hspace{0.1in}
   \subfigure[Quadrature points (run 4) ]{ 
    \label{test2:aqr} 
    \includegraphics[width=1.3in]{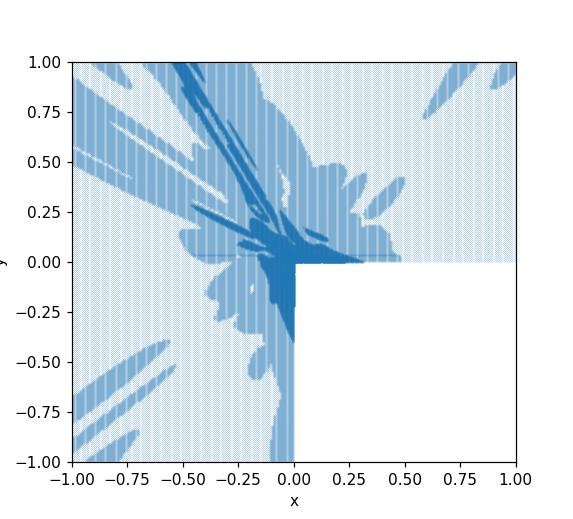}}  
     \hspace{0.1in}
  \subfigure[Marked elements at run 2]{ 
    \label{test2:mark1} 
    \includegraphics[width=1.5in]{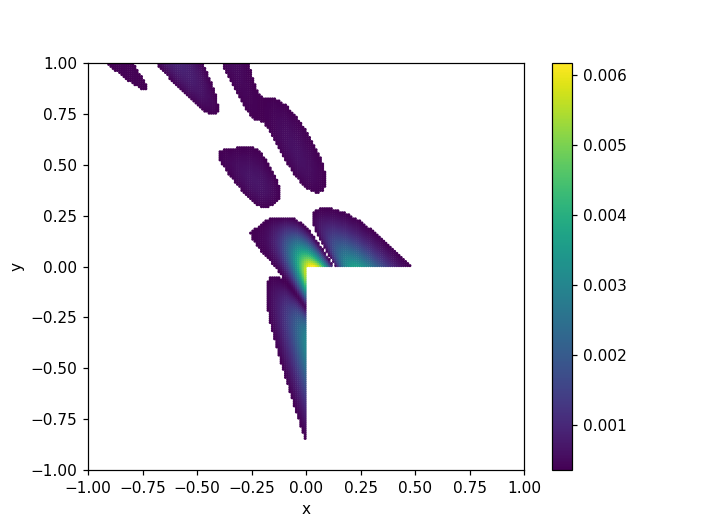}}
     \hspace{0.1in}
   \subfigure[Marked elements at run 3]{ 
    \label{test2:mark2} 
    \includegraphics[width=1.5in]{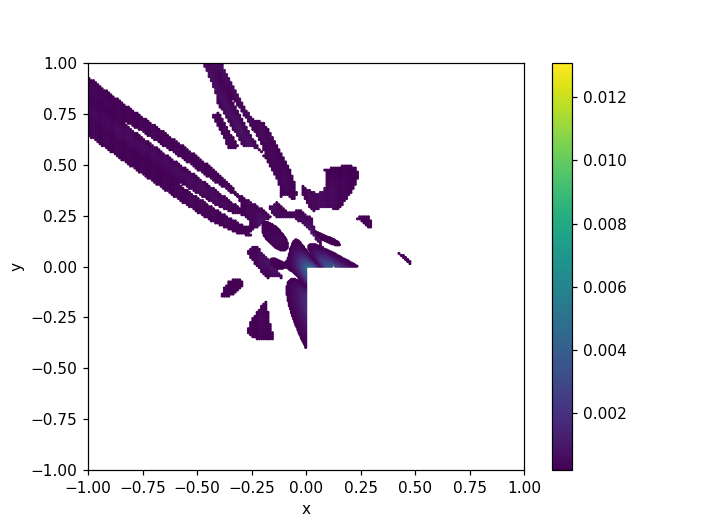}}
     \hspace{0.1in}
  \subfigure[Marked elements at run 4]{ 
    \label{test2:mark3} 
    \includegraphics[width=1.5in]{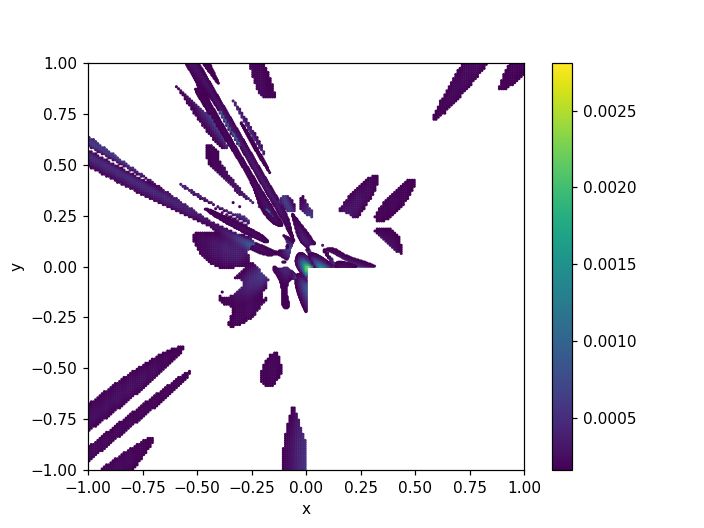}} 

  \caption{Test case II numerical solution using adaptive deep Ritz (2-48-48-48-2). (a-d) component-wise numerical solution $\vu$ and $\sigma$; (f)final quadrature points obtained through AQR; (g-i) marked element during the four-run AQR process. }
  \label{Lshape_results} 
\end{figure}

Due to a stress singularity at the corner point $(0,0)$, the direct LS physics-driven methods do not apply to this type of problems. Using the deep Ritz method, we test the performance of a four layer DNN structure (2-48-48-48-2, 4898 parameters) with fixed uniform quadrature and adaptive quadrature refinement (AQR). As shown in Table. \ref{test2_table}, with AQR generated non-uniformly distributed quadrature points, adaptive deep Ritz approximates the solution using less data (quadatrure points) compared with the uniform fixed quadrature method. During the iterative process, the marked regions with larger residuals who need refinement are depicted in Figure. \ref{test2:mark1} -  \ref{test2:mark3} and the generated non-uniform quadrature points are plotted in Figure. \ref{test2:aqr}. The stress singularity is well captured by the local error indicators and correspondingly, more quadrature points are distributed over the re-entrant corner region. The final numerical solutions obtained are plotted in Figure. \ref{test2:u1} -  \ref{test2:s22}. This experiment shows the validity of the proposed local error estimator for the adaptive quadrature scheme.

\begin{table}[htbp]
\centering
\caption{Numerical results of adaptive deep Ritz for test case II using a DNN structure 2-48-48-48-2.}
\vspace{5pt}
\begin{tabular}{|c|c|c|c|c|c|}
\hline
{\bfseries \makecell{Quadrature\\ method}}  & \multicolumn{2}{c|}{ $\#\cT$}
& $\dfrac{\|\vu-\vu_{_N}\|_a}{\|\vu\|_a}$ &  $\dfrac{\|\sigma-\sigma_{_N}\|}{\|\sigma\|}$ &  $\dfrac{\|\vu-\vu_{_N}\|}{\|\vu\|}$ \\
\hline
{\makecell{uniform \\fixed}} &\multicolumn{2}{c|}{$120,000$} &10.99\%  & 10.20\%  & 1.83\% \\
\hline
\multirow{4}{*}{\makecell{non-uniform \\AQR}} & run 1 & 30,000 &23.38\%  &21.58\% &3.31\%  \\
\cline{2-6}
& run 2 &42,897 &14.44\%    &13.45\%  & 1.78\% \\
\cline{2-6}
& run 3 &60,306  &11.32\%  & 10.48\%  &1.69\%   \\
\cline{2-6}
& run 4 &90,237 &10.71\%  &9.86\%   & 1.68\%  \\
\hline
\multicolumn{6}{l}{\small *Training details:  $\gamma_{_D}= 1$ ; } \\
\multicolumn{6}{l}{\small for the uniform quadrature and the AQR  run 1: trained with $200,000$ iterations with } \\
\multicolumn{6}{l}{\small learning rate starts from 0.01 and decays $90\%$ every 50,000 iterations;} \\
\multicolumn{6}{l}{\small for non-uniform AQR run 2-4, trained with 100,000 iterations using fixed learning rate 1e-5.} \\
\end{tabular}
\label{test2_table}
\end{table}

\subsection{Test case III: A quadratic membrane under tension}
\begin{figure}[htb!]
  \centering 
  \subfigure[The test problem]{ 
    \label{test3:problem} 
    \includegraphics[width=1.3in]{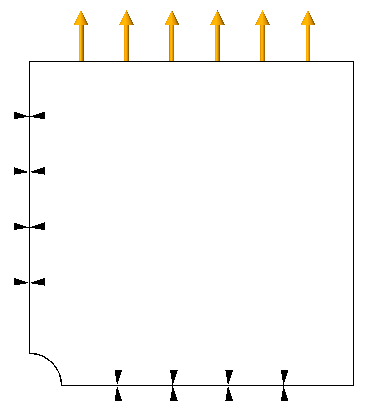}} 
    \hspace{0.2in}
    \subfigure[FEA $\tens{\sigma}_{yy}$ ]{ 
    \label{test3:fea_syy} 
    \includegraphics[width=1.5in]
    {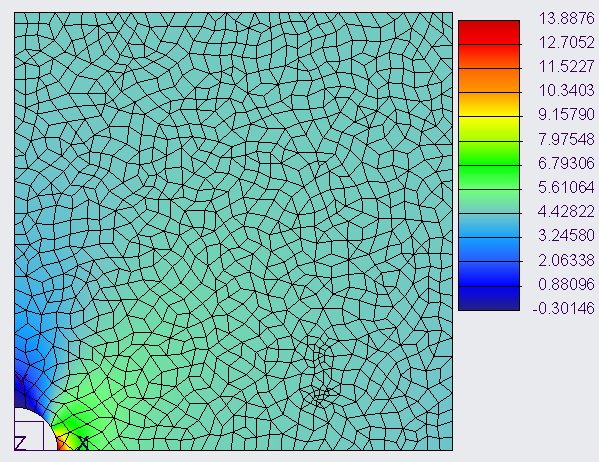}}
    \hspace{0.15in}
    \subfigure[FEA $\vu_{y}$ ]{ 
    \label{test3:fea_uy} 
    \includegraphics[width=1.5in]{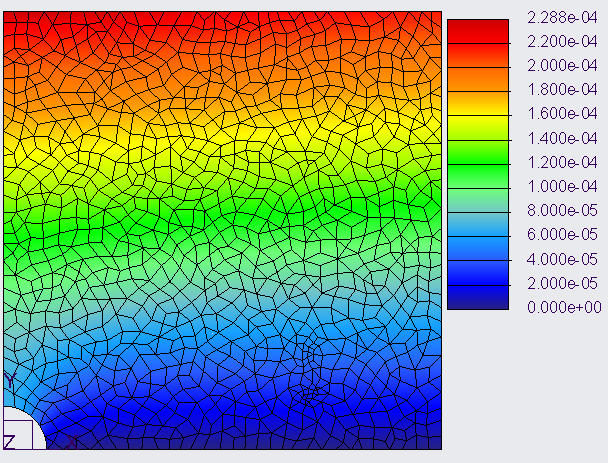}}
    \hspace{0.15in}
  \subfigure[Last run AQR points]{
    \label{test3:aqr} 
    \includegraphics[width=1.3in]{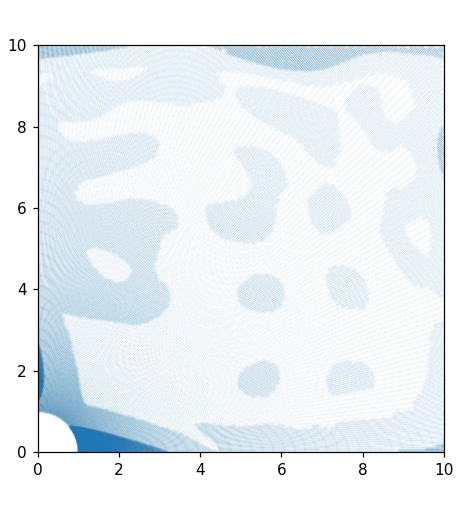}}
        \hspace{0.15in}
   \subfigure[Ritz $\tens{\sigma}_{yy}$ ]{ 
    \label{test3:ritz_sol} 
    \includegraphics[height=1.35in]{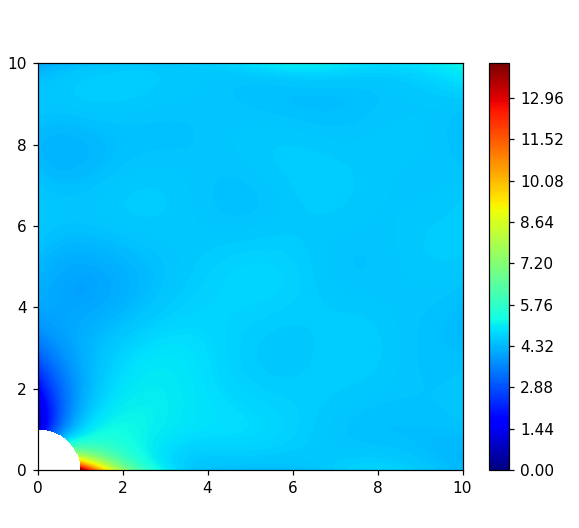}}
        \hspace{0.15in}
   \subfigure[Ritz $\vu_{y}$]{ 
    \label{test3:ritz:uy} 
    \includegraphics[width=1.6in]{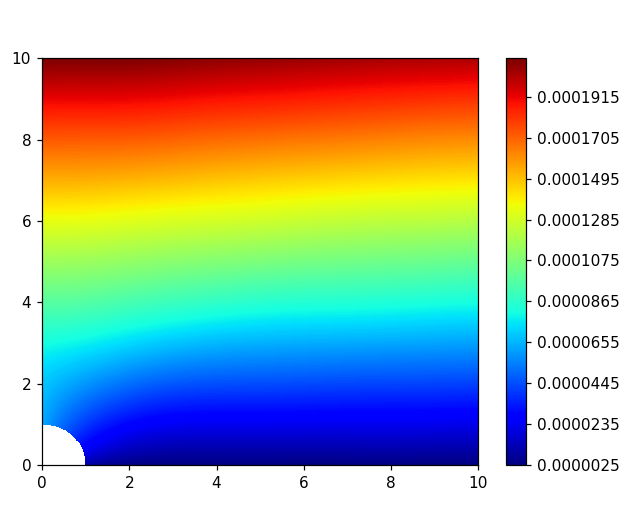}}
    
  \caption{Test case III Numerical solution using adatpive FEA p-element and deep Ritz with adaptive quadrature. 
    (Rtiz solution is obtained with a DNN structure: 2-64-64-64-2 of 8578 DoF. FEA solution is obtained using p-element with 14786 DoF.)}
  \label{Result_3_Ritz} 
\end{figure}

The third test problem is given by a quadratic membrane of elastic isotropic material with a circular hole in the center. Traction forces act on the upper and lower edges of the strip, body forces are ignored. Because of the
symmetry of the problem, it suffices to compute only a fourth of the total geometry. The computational domain is then given by
$$
\Omega = \{\vx \in \RR^2: 0< x <10, 0<y <10, x^2+y^2>1\}.
$$
The boundary condition on the top edge of the computation domain ($\G_1: \{y=10,  0< x <10\}$) are set to $\tens{\sigma}{\vn} = (0,4.5)^T$, the boundary condition on the bottom ($\G_2: \{y=0,  1< x <10\}$) are set to $(\sigma_{xx}, \sigma_{xy})\cdot {\vn} = 0, u_y =0$ (symmetry condition), and finally, the boundary condition on the left ($\G_3: \{x=0,  1< y <10\}$) are given by $(\sigma_{yx}, \sigma_{yy})\cdot {\vn} = 0$, and $u_x =0$ (symmetry condition). The material parameters are $E= 206900$ for Young's modulus and and $\nu=0.29$ for Poisson's ratio.
The challenge of this test is the stress concentration located around point (0,1) due to the presence of the small hole. Since there is no analytic solution, a reference solutions is obtained using finite element analysis (FEA) with an adaptive mesh refinement (adaptive p-element refined with highest polynomial order of 7). The reference solution is given in Figures \ref{test3:fea_syy} and \ref{test3:fea_uy}, where we use two key measures:  $\max{\tens{\sigma}_{yy}} = 13.8876$,  an $\max{\|\vu\|}= $2.288e-4 as references (a similar reference solution is also provided in \cite{Cai04} using adaptive h-elements). 

 \begin{table}[htbp]
\centering
\caption{Using Deep Ritz with adaptive quadrature to solve test case III}
\vspace{5pt}
\begin{tabular}{|l|c|c|c|c|}
\hline

\textbf{Iteration}  & \textbf{No. of Quad} & \textbf{max $\|\vu\|$ }  &\textbf{max} $\tens{\sigma}_{yy}$ &$\dfrac{\sum_{T\in \cT}\eta_{_T}}{\#\cT}$ \\ \hline
1 &  28,392 & 2.1327e-4  & 7.9155 & 9.5942e-3 \\\hline 
2 & 40,878 & 2.1566e-4  & 10.0728 & 6.7179e-4 \\\hline
3 & 60,434 & 2.1623e-4  &13.2632  & 4.4317e-4  \\\hline
4 & 88,574 & 2.2672e-4 & 13.8912 & 2.9309e-4  \\\hline
\multicolumn{5}{l}{\small *Training details: $\gamma_{_D}= 1$;} \\
\multicolumn{5}{l}{\small run 1 is trained with $200,000$ iterations, with learning rate starts from} \\
\multicolumn{5}{l}{\small 0.01 and decays $90\%$ every 50,000 iterations;} \\
\multicolumn{5}{l}{\small run 2-4: 50,000 iterations each with fixed learning rate 1e-5.} \\
\end{tabular}
\centering
\label{test3_table1}
\end{table}



We tested the adaptive deep Ritz with a four-layer structure (2-64-64-64-2, 8578 parameters). Initially, 28,392 uniform quadrature points using polar coordinates are used for getting the first iteration that captures well for the displacement field $\vu$.  However the stress concentration factor is not accurate with uniformly distributed quadrature points. After four iterations of AQR using the average marking strategy ($\gamma_2 = 1.5)$, a total of 88,574 quadrature points are generated adaptively as shown in Figure \ref{test3:aqr} and with this set of non-uniform quadrature points, the method is capable of obtaining a similar stress concentration factor compared with the adaptive FEA solution, see the numerical results depicted in Figure. \ref{test1ritz:a}-\ref{test1ritz:c} and Table. \ref{test3_table1}. As expected, more quadrature points are distributed adaptively near the small hole region during the AQR processes such that the stress concentration can be simulated accurately.

\begin{figure}[ht]
\centering
\includegraphics[ height=3.8in]{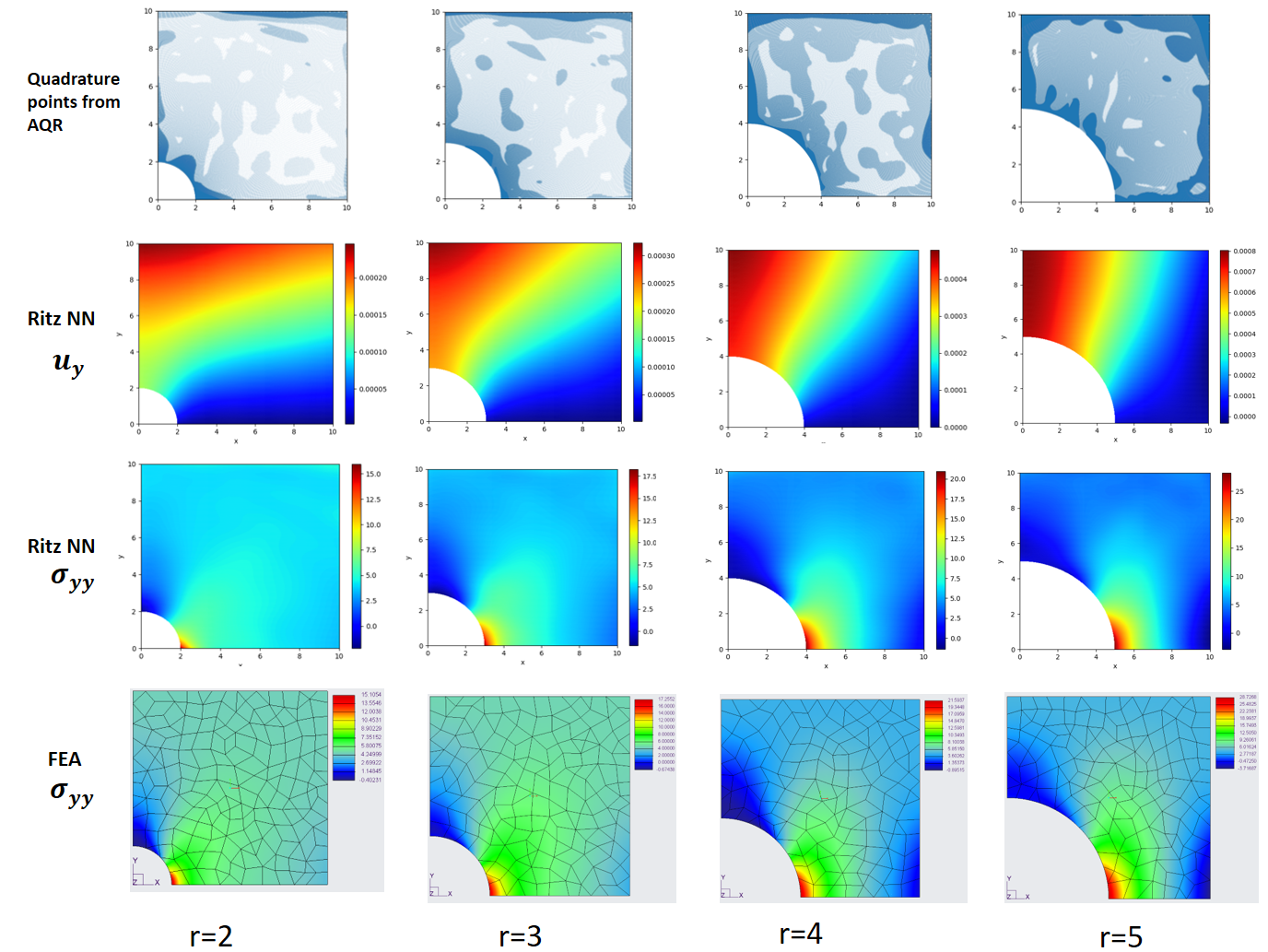}
\caption{Sensitivity study of the test case III with varying hole radii. First row shows the generated adaptive quadrature points at the listed hole radius steps; second and third rows show the numerical solution using adaptive deep Ritz, and the last row give the reference adaptive FEA solutions.}
\label{test3_sensitivity}
\end{figure}

Furthermore, we explored the potential of using transfer learning to solve a family of stress problems. To this end, we conducted a sensitivity study by varying the center hole from $r=1$ to $r=5$, taking a step size of $0.5$. Using the trained DNN model in the previous experiment (with an initial hole size of $r=1$) as the starting point, we stepped through the various hole sizes by training a DNN using weight transferring. Our assumption was that the DNN trained from the previous step forms a good initial for the next parametric step, and therefore the adaptive deep Ritz would converge faster if applied to a family of similar problems.  In our test, we verified this assumption and it took significantly fewer number of iterations (in this example, 50,000 iterations and two-run AQR for each hole size step) to converge to the results.  This is compared to the total 350,000 iterations used for the result obtained in Table . \ref{test3_table1},  demonstrating that transfer learning saves significant time.  Our results at each step also align with the numerical solutions evaluated through FEA adaptive p-refinement method. Figure. \ref{test3_sensitivity} lists the displacement component $\vu_y$ and the associate stress tensor component $\tens{\sigma}_{yy}$ obtained from adptive deep Ritz transfer learning and adaptive FEA. 

 \subsection{Discussion}
 \textbf{Numerical Differentiation (ND) or Automatic Differentiation (AD).}
For each quadrature point $\vx_{_T}$, the evaluations of  $\tens{\varepsilon}\big(\vv(\vx_{_T})\big)$ and $\nabla\!\cdot\!\vv(\vx_{_T})$ in the energy functional are based on the first order partial derivatives that may be calculated through either automatic differentiation (AD) or numerical differentiation (ND). AD is a common choice for some physics-driven methods such as the PINN and its variants. When a DNN function has the first order derivatives at all sampling points, AD eliminates numerical error caused by ND. However, when using AD, there are some difficulties in training PINN as noticed in \cite{CHIU2022}; moreover, a discrete differential operator combining AD and ND was proposed.
In our experiments, pathologies in training deep Ritz with AD were also observed. Nevertheless, our experiments show numerically that ND with quadrature-based numerical integration is robust in training/solving the deep Ritz.

\begin{figure}[htb!]
  \centering 
  \subfigure[Energy Functional of different methods]{ 
    \label{test1:energy} 
    \includegraphics[width=2.8in]{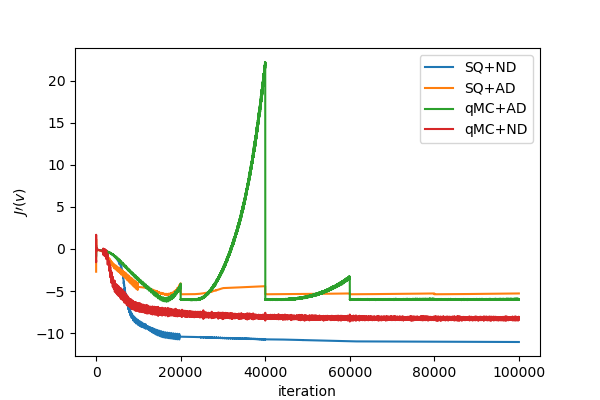}} 
    \hspace{0.1in}
  \subfigure[$\|\vu\|$ (SQ+ND)]{ 
    \label{test1:sqnd} 
    \includegraphics[width=1.6in]{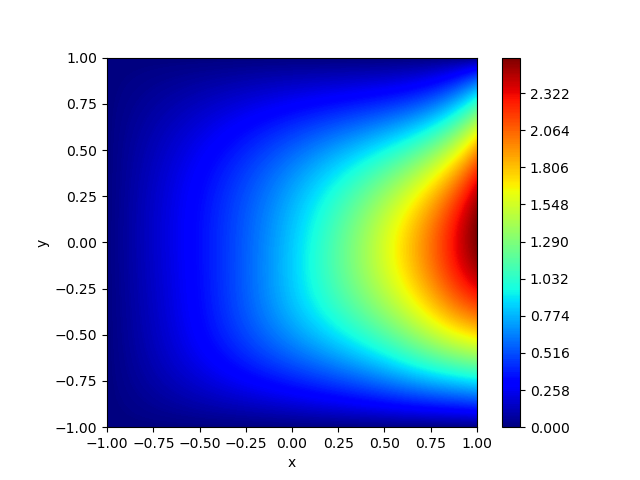}}
    \hspace{0.1in}
    \subfigure[$\|\vu\|$ (SQ+AD)]{ 
    \label{test1:sqad} 
    \includegraphics[width=1.6in]{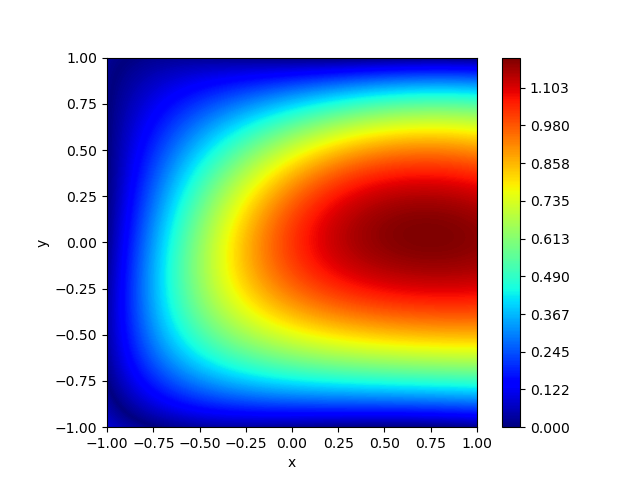}}
    \hspace{0.1in}
  \subfigure[$\|\vu\|$ (qMC+ND)]{
    \label{test1:mcnd} 
    \includegraphics[width=1.6in]{Fig/test1_MC_ND_U.png}}  
    \hspace{0.1in}
  \subfigure[$\|\vu\|$ (qMC+AD)]{
    \label{test1:mcad} 
    \includegraphics[width=1.6in]{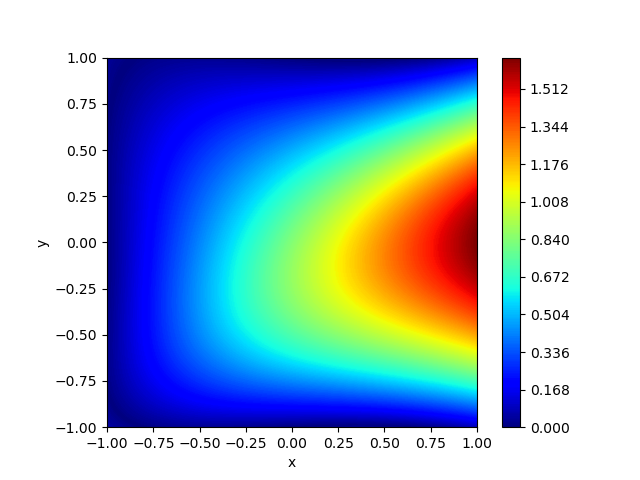}}
    
  \caption{Comparison study for test case I using numerical differentiation v.s. automatic differentiation, and standard quadrature-based integration v.s. quasi-Monte Carlo-based integration. All results are trained using deep Ritz DNN structure (2-32-32-2).}
  \label{Result_1_NDAD} 
\end{figure}

Here, we use the  test case I as an example and list the results of using combinations of different integration and differentiation methods: (1) standard quadrature-based integration plus AD (SQ+AD),  (2) standard quadrature-based integration plus ND (SQ+ND), (3) quasi-Monte Carlo-based integration plus AD (qMC+AD), and (4) quasi-Monte Carlo-based integration plus ND (qMC+ND). As shown in Figure. \ref{test1:energy}, only SQ+ND converges to the real potential energy of the tested problem $J^{\star}(\vu)\approx -11.2$ within 100,000 iterations, the other three methods are suspiciously trapped in local minima and the obtained solutions are non-physical as shown in Figure. \ref{test1:sqad}-\ref{test1:mcad}. 


\textcolor{black}{\textbf{Non-convex optimization}  The exceptional approximation power of neural networks allows them to effectively represent a wide range of solutions, including those with irregular geometries, discontinuities, or singularities that can pose challenges for traditional finite element methods.  However, this type of approximating functions also introduces a computationally demanding optimization problem. NN-based algorithms for solving PDEs typically involve a high-dimensional and non-convex optimization for which the first order stochastic gradient descent-based methods are the most widely employed solvers. Currently, NN-based methods are not competitive with well-established mesh-based methods due to the considerable computational cost of the algebraic solvers, even though NN-based methods may require fewer degrees of freedom than their mesh-based counterparts. Developing fast solvers is an open and challenging problem and requires lots of efforts from numerical
analysts. }

\section{Conclusion}
In this paper, linear elasticity problems are formulated under the Ritz framework and are discretized using DNN functions. To enforce the essential boundary condition, the energy functional is modified with an extra penalization term using $H^{1/2}$ norm. It is shown that within the function class, the minimization of the modified energy functional yields the best approximation with respect to the modified energy norm. To calculate the modified energy functional accurately and efficiently, adaptive quadrature refinement equipped with a local residual-based error indicator was proposed and tested, and its effectiveness and efficiency in improving numerical simulation was demonstrated.


There are still numerous unresolved issues that require further research. In this study, we make the assumption that the DNN is sufficiently large to approximate the solution. However, selecting an appropriate network structure for different problems and establishing a suitable initial DNN model are still open questions. While we conducted a basic sensitivity analysis to demonstrate the potential of DNNs for parametric PDEs, exploring efficient methods to extend the transfer learning strategy to general design space exploration and even solving topology optimization problems requires further investigation.   



\section*{Acknowledgments}

This work was supported in part by the National Science Foundation under grant DMS-2110571 and the Future of Work at the Human Technology Frontier (FW-HTF) 1839971. We also acknowledge the Feddersen Distinguished Professorship Funds. Any opinions, findings, and conclusions expressed in this material are those of the authors and do not necessarily reflect the views of the funding agency. We thank the support of DARPA project on symbiotic design and Stanford Research International (SRI) for partial support of the project.

\bibliographystyle{elsarticle-num}
\bibliography{reference}

\begin{thebibliography}{10}
\expandafter\ifx\csname url\endcsname\relax
  \def\url#1{\texttt{#1}}\fi
\expandafter\ifx\csname urlprefix\endcsname\relax\def\urlprefix{URL }\fi
\expandafter\ifx\csname href\endcsname\relax
  \def\href#1#2{#2} \def\path#1{#1}\fi

\bibitem{Liang18}
L.~Liang, M.~Liu, C.~Martin, W.~Sun, A deep learning approach to estimate
  stress distribution: a fast and accurate surrogate of finite-element
  analysis, Journal of The Royal Society Interface 15~(138) (2018) 20170844.

\bibitem{Gao20}
W.~Gao, X.~Lu, Y.~Peng, L.~Wu, A deep learning approach replacing the finite
  difference method for in situ stress prediction, IEEE Access 8 (2020)
  44063--44074.

\bibitem{Wang21}
S.~Wang, Y.~Teng, P.~Perdikaris, Understanding and mitigating gradient flow
  pathologies in physics-informed neural networks, SIAM Journal on Scientific
  Computing 43~(5) (2021) A3055--A3081.

\bibitem{Badarinath21}
P.~Vurtur~Badarinath, M.~Chierichetti, F.~Davoudi~Kakhki, A machine learning
  approach as a surrogate for a finite element analysis: Status of research and
  application to one dimensional systems, Sensors 21~(5) (2021).

\bibitem{iakovlev2021learning}
V.~Iakovlev, M.~Heinonen, H.~L{\"a}hdesm{\"a}ki, Learning continuous-time
  {PDE}s from sparse data with graph neural networks, in: International
  Conference on Learning Representations, 2021.

\bibitem{li2021fourier}
Z.~Li, N.~B. Kovachki, K.~Azizzadenesheli, B.~Liu, K.~Bhattacharya, A.~Stuart,
  A.~Anandkumar, Fourier neural operator for parametric partial differential
  equations, in: International Conference on Learning Representations, 2021.

\bibitem{raissi2017hidden}
M.~Raissi, G.~E. Karniadakis, Hidden physics models: Machine learning of
  nonlinear partial differential equations, Journal of Computational Physics
  (2017).

\bibitem{Long19}
Z.~Long, Y.~Lu, B.~Dong, {PDE}-net 2.0: Learning pdes from data with a
  numeric-symbolic hybrid deep network, J. Comput. Phys. 399 (2019).

\bibitem{Khoo19}
Y.~Khoo, L.~Ying, Switchnet: A neural network model for forward and inverse
  scattering problems, SIAM Journal on Scientific Computing 41~(5) (2019)
  A3182--A3201.

\bibitem{Weinan17}
W.~E, J.~Han, A.~Jentzen, Deep learning-based numerical methods for
  high-dimensional parabolic partial differential equations and backward
  stochastic differential equations, Communications in Mathematics and
  Statistics 5~(4) (2017) 349--380.

\bibitem{Sirignano18}
J.~Sirignano, K.~Spiliopoulos, {DGM}: A deep learning algorithm for solving
  partial differential equations, Journal of Computational Physics 375 (2018)
  1139--1364.

\bibitem{Berg18}
J.~Berg, K.~Nystrom, A unified deep artificial neural network approach to
  partial differential equations in complex geometries, Neurocomputing 317
  (2018) 28--41.

\bibitem{Weinan18}
W.~E, B.~Yu, The deep {R}itz method: A deep learning-based numerical algorithm
  for solving variational problems, Communications in Mathematics and
  Statistics 6~(1) (2018) 1--12.

\bibitem{Karniadakis19}
M.~Raissi, P.~Perdikaris, G.~E. Karniadakis, Physics-informed neural networks:
  A deep learning framework for solving forward and inverse problems involving
  nonlinear partial differential equations, Journal of Computational Physics
  378 (2019) 686–707.

\bibitem{CAI2020}
Z.~Cai, J.~Chen, M.~Liu, X.~Liu, Deep least-squares methods: An unsupervised
  learning-based numerical method for solving elliptic {PDE}s, Journal of
  Computational Physics 420 (2020) 109707.

\bibitem{Cai2021linear}
Z.~Cai, J.~Chen, M.~Liu, Least-squares {R}e{LU} neural network {(LSNN)} method
  for linear advection-reaction equation, Journal of Computational Physics 443
  (2021) 110514.

\bibitem{LI-wei21}
W.~Li, M.~Z. Bazant, J.~Zhu, A physics-guided neural network framework for
  elastic plates: Comparison of governing equations-based and energy-based
  approaches, Computer Methods in Applied Mechanics and Engineering 383 (2021)
  113933.

\bibitem{LiuCai1}
M.~Liu, Z.~Cai, J.~Chen, Adaptive two-layer {ReLU} neural network: I. best
  least-squares approximation, Computers. Math. Appl. 113 (2022) 34--44.

\bibitem{CAI2022}
Z.~Cai, J.~Chen, M.~Liu, Least-squares {R}e{LU} neural network {(LSNN)} method
  for scalar nonlinear hyperbolic conservation law, Applied Numerical
  Mathematics 174 (2022) 163--176.

\bibitem{Cai2023}
Z.~Cai, J.~Chen, M.~Liu, Least-squares {ReLU} neural network {(LSNN)} method
  for scalar nonlinear hyperbolic conservation laws: discrete divergence
  operator, J. Comput. Appl. Math. (2023) to appear.

\bibitem{Xu2020}
J.~Xu, The finite neuron method and convergence analysis, Communications in
  Computational Physics 28 (2020) 1707--1745.

\bibitem{LiuCai2}
M.~Liu, Z.~Cai, Adaptive two-layer {R}e{LU} neural network: {II}. {Ritz}
  approximation to elliptic {PDEs}, Computers. Math. Appl. 113 (2022) 103--116.

\bibitem{bg5}
P.~B. Bochev, M.~D. Gunzburger, Least-Squares Finite Element Methods, Vol. 166,
  Applied Mathematics Sciences, Springer, 2009.

\bibitem{Cai04}
Z.~Cai, G.~Starke, Least-squares methods for linear elasticity, SIAM Journal on
  Numerical Analysis 42~(2) (2004) 826--842.

\bibitem{BeCaPa:19}
F.~Bertrand, Z.~Cai, E.-Y. Park, Least-squares methods for elasticity and
  stokes equations with weakly imposed symmetry, Comput. Methods Appl. Math
  19~(3) (2019) 415--430.

\bibitem{HAGHIGHAT2021}
E.~Haghighat, M.~Raissi, A.~Moure, H.~Gomez, R.~Juanes, A physics-informed deep
  learning framework for inversion and surrogate modeling in solid mechanics,
  Computer Methods in Applied Mechanics and Engineering 379 (2021) 113741.

\bibitem{Nitsche71}
J.~Nitsche, Über ein variationsprinzip zur lösung von dirichlet-problemen bei
  verwendung von teilräumen, die keinen randbedingungen unterworfen sind, Abh.
  Math. Sem. Univ. 36~(1) (1971) 9--15.

\bibitem{Liao21}
Y.~Liao, P.~Ming, Deep {N}itsche method: Deep {R}itz method with essential
  boundary conditions, Communications in Computational Physics 29~(5) (2021)
  1365--1384.

\bibitem{URIARTE23}
C.~Uriarte, D.~Pardo, I.~Muga, J.~Muñoz-Matute, A deep double {R}itz method
  (d2rm) for solving partial differential equations using neural networks,
  Computer Methods in Applied Mechanics and Engineering 405 (2023) 115892.

\bibitem{RIVERA22}
J.~A. Rivera, J.~M. Taylor, Ángel J.~Omella, D.~Pardo, On quadrature rules for
  solving partial differential equations using neural networks, Computer
  Methods in Applied Mechanics and Engineering 393 (2022) 114710.

\bibitem{pinkus1999}
A.~Pinkus, Approximation theory of the mlp model in nueral networks, Acta
  Numerica 15 (1999) 143--195.

\bibitem{Cybenko1989}
G.~Cybenko, Approximation by superpositions of a sigmoidal function,
  Mathematics of Control, Signals, and Systems 2 (1989) 303--314.

\bibitem{HornikS1989}
K.~Hornik, M.~Stinchcombe, H.~White, Multilayer feedforward networks are
  universal approximators, Neural Networks 2 (1989) 359--366.

\bibitem{Schumaker}
L.~Schumaker, Spline Functions: Basic Theory, John Wiley, New York, 1981.

\bibitem{Ciarlet78}
P.~G. Ciarlet, The finite element method for elliptic problems, Society for
  Industrial and Applied Mathematics, 1978.

\bibitem{kingma2015}
D.~P. Kingma, J.~Ba, Adam: A method for stochastic optimization, in:
  International Conference on Representation Learning, San Diego, 2015.

\bibitem{Harper19}
G.~Harper, J.~Liu, S.~Tavener, B.~Zheng, Lowest-order weak {G}alerkin finite
  element methods for linear elasticity on rectangular and brick meshes, J.
  Sci. Comput. 78~(3) (2019) 1917–1941.

\bibitem{Alberty02}
J.~Alberty, C.~Carstensen, S.~A. Funken, R.~Klose, Matlab implementation of the
  finite element method in elasticity, Computing 69~(3) (2002) 239--263.

\bibitem{CHIU2022}
P.-H. Chiu, J.~C. Wong, C.~Ooi, M.~H. Dao, Y.-S. Ong, Can-pinn: A fast
  physics-informed neural network based on coupled-automatic–numerical
  differentiation method, Computer Methods in Applied Mechanics and Engineering
  395 (2022) 114909.

\end{thebibliography}

\end{document}